\documentstyle[a4,amstex]{article}
\newcommand{\real}{{\Bbb R}}
\newcommand{\cplx}{{\Bbb C}}
\newcommand{\zint}{{\Bbb Z}}
\newcommand{\cplxn}{\cplx^{n}}
\newcommand{\ep}{\epsilon}
\newcommand{\Hhat}{\widehat{H}_N(q)}

\newcommand{\sll}{\widehat{ \frak{s}\frak{l}}_n}
\newcommand{\sln}{ {{\frak s} {\frak l}}_n}
\newcommand{\slt}{\widehat{ \frak{s}\frak{l}}_2}
\newcommand{\tl}{ \frak{s}\frak{l}_2}
\newcommand{\lb}{\boldsymbol{ \lambda}}
\newcommand{\mb}{\boldsymbol{ \mu}}
\newcommand{\mm}{{\bold m}}
\newcommand{\nn}{{\bold n}}
\newcommand{\zz}{{\bold z}}
\newcommand{\ee}{{\bold e}}
\newcommand{\ro}{\rho^M_{r+1,r}}
\newcommand{\rk}{\rho^{M,k}_{r+1,r}}
\newcommand{\MC}{{\cal M}}
\newcommand{\EC}{{\cal E}}
\newcommand{\FC}{{\cal F}}

\renewcommand{\l}{\lambda}
\newcommand{\s}{\sigma}

\newcommand{\p}{\partial}
\newcommand{\ba}{\begin{array}}
\newcommand{\ea}{\end{array}}
\newcommand{\beq}{\begin{equation}}
\newcommand{\eeq}{\end{equation}}
\newcommand{\bqa}{\begin{eqnarray}}
\newcommand{\eqa}{\end{eqnarray}}
\newcommand{\bqas}{\begin{eqnarray*}}
\newcommand{\eqas}{\end{eqnarray*}}
\newcommand{\halmos}{\rule{5pt}{5pt}}

\numberwithin{equation}{section}
\newtheorem{df}{\bf Definition}
\newtheorem{prop}{\bf Proposition}

\newtheorem{lemma}{\bf Lemma}
\newtheorem{cor}{\bf Corollary}
\newenvironment{rmk}{\noindent{\bf Remark}\hskip 5pt}{}

\newenvironment{pf}{\noindent {\em {\normalsize P}\normalsize{roof.}}  
\normalsize\hskip 5pt}{\hfill\halmos}
\newenvironment{ppf}[1]{\noindent {\bf Proof of the Proposition #1}  
\normalsize\hskip 5pt}{\hfill\halmos}

\begin{document}

\begin{titlepage}
\pagestyle{empty}
\begin{center}
\mbox{} \\
\vspace{2.7cm}
\begin{Large}
{\bf Level-0 action of $U_q(\sll)$ on the $q$-deformed Fock spaces } 
\end{Large} \\
\vspace{1.5cm}
\large{ Kouichi Takemura }\footnote{e-mail: takemura@@kurims.kyoto-u.ac.jp} 
\large{ and Denis Uglov} \footnote{e-mail: duglov@@kurims.kyoto-u.ac.jp}  
\\  Research Institute for Mathematical Sciences\\
Kyoto University, Kyoto 606, Japan  \\ 
\vspace{1cm}
July 1996  \\

\vspace{2cm}

\begin{abstract}
On the level-1 Fock space modules of the algebra $U_q(\sll)$ we define a level-0 action $U_0$ of the $U_q(\sll)$, and an action of an abelian algebra of conserved Hamiltonians  commuting with the $U_0$. An irreducible decomposition of the Fock space with respect to the level-0 action is derived by constructing a base of the Fock space in terms of the Non-symmetric Macdonald Polynomials. 

\end{abstract}
\end{center}
\end{titlepage}

\voffset=-1.5cm
\hoffset=-1.5cm

\section{Introduction}

Recently several intriguing links have been uncovered between Solvable Models with Long-Range Interaction and representations of affine Lie algebras. 

The earliest example is the identification of the Field Theory limit in the Haldane-Shastry spin chain with the level-1 $su(2)$ WZW Conformal Field Theory made in \cite{HHTBP}. This led to discoveries such as the Yangian action on the integrable level-1 modules of the affine Lie algebra $\sll$, the spinon bases and the fermionic character formulas \cite{BPS,BBL}. 

Another line of research connected the Calogero-Sutherland Model and its $q$-deformed analog -- the trigonometric Ruijsenaars Model -- with the $W$-algebras and their $q$-deformations \cite{Awata1}.     

In the work of Bernard {\em et al.} \cite{BG} both the Haldane-Shastry and the Calogero-Sutherland Models were understood to be special cases of the more  general $\sln$-invariant Dynamical or Spin Calogero-Sutherland Model. In this light it becomes plausible that a connection between the Long-Range Models and CFT must exist  on  the level of the Spin Calogero-Sutherland Model as well, so that the two lines of results mentioned above can be seen as parts of a more general structure.

Two attempts to understand the Spin Calogero-Sutherland Model in the  Field Theory limit were made by different means  in $\cite{Awata2}$ and in $\cite{U}$. In $\cite{U}$ this limit was constructed for the $\tl$-invariant case by using the semi-infinite wedge realization of the Fock space module of the affine Lie algebra $\slt$ \cite{KacRaina}. The Fock space was interpreted as the space of states of the Model in the Field Theory limit, it  was shown to admit a Yangian action and the decomposition of this space with respect to this action was derived.

The aim of the present paper is to give a $q$-deformation of the construction in $\cite{U}$ in the $\sln$-invariant situation. Let us briefly go through the main features of our work. We start with the $q$-deformation of the $N$-particle fermionic Spin Calogero-Sutherland Model which was proposed in \cite{BG} (see also \cite{Konno1} ). The ingredients that define the Model are: a family of mutually commuting operators -- conserved Hamiltonians $ h^{(N)}_l $,  $(l=\pm1,\pm2,\dots \;)$ and the level-0 action of the $q$-deformed affine algebra $U_q(\sll)$ which we denote by $U_0^{(N)}$. This action commutes with the operators $ h^{(N)}_l $ and hence has a meaning of a non-abelian symmetry algebra of the Model. The space of states in the Model is identified with the finite $N$-fold $q$-wedge product which was recently introduced in \cite{KMS}, it is isomorphic to the space of $q$-fermionic states proposed in \cite{BG}. 

Our main problem can be formulated as that of taking the number of particles $N$ in the Model to infinity so that the commuting Hamiltonians and the symmetry algebra remain well-defined. To solve this problem we utilize a certain projective limit of the finite $q$-wedge product -- the semi-infinite $q$-wedge product of \cite{KMS}.  

The  semi-infinite $q$-wedge product is shown in \cite{KMS} to be isomorphic to the level-1 $q$-deformed Fock space module of $U_q(\sll)$ discovered in \cite{H}. We are able to demonstrate that the Hamiltonians and the commuting with them level-0 action of $U_q(\sll)$ carry over from the finite into the semi-infinite $q$-wedge product and therefore define the Model in the Fock space. By analogy with the rational case considered in \cite{U} we interpret this as a Field Theory limit of the $q$-deformed Spin Calogero-Sutherland Model and of the associated symmetry algebra. 

The second problem which we consider is to derive the decomposition of the Fock space with respect to the level-0 action, and diagonalize the commuting Hamiltonians. Towards this end we construct a suitable base  of the Fock space in terms of the Non-symmetric Macdonald Polynomials \cite{Macdonald,Cherednik}. One of the features of the level-0 decomposition of the Fock space is its relative simplicity  -- at generic values of parameters in the Model each irreducible component is isomorphic to a tensor product of {\em fundamental} representations of $U_q(\sll)$ in terminology of \cite{CP} -- this makes it different from the level-0 decomposition of the irreducible level-1 $U_q(\slt)$ modules constructed in \cite{J} where an irreducible component is in general a subquotient of a tensor product of the fundamental representations.

The arrangement of the paper is as follows. The sec.2 is has an introductory character, here we give the relevant background information on the affine Hecke algebra and $q$-wedge products. In sec. 3 the space of finite wedges is decomposed with respect to the  $U_q(\sll)$-action $U_0^{(N)}$. In sec. 4 we define the level-0 action and the family of commuting Hamiltonians in the Fock space. Sec. 5 contains results on the decomposition of the Fock space with respect to this action. 

\noindent {\large {\bf Acknowledgments} } We would like to thank Professors T.Miwa and M.Kashiwara for discussions and support.

\section{ Preliminaries}
In this introductory section we set up our notations and collect several known definitions and results to be used starting from sec. 3. For details one may consult the works \cite{KMS,BG,Cherednik,Macdonald}. 
\\ \mbox{} \\

\subsection{ Affine Hecke algebra}  The affine Hecke algebra $\Hhat$ is an associative algebra generated by elements $ T_i $  ( $ i=1,\dots, N-1 $) and $ y_j^{\pm 1} $  $( j=1,\dots, N) $. These elements satisfy the following relations:      
\bqa
T_i^2 & = & (q - q^{-1}) T_i + I , \\
T_i T_{i+1} T_i & = & T_{i+1} T_i T_{i+1} , \\
T_i T_j & = & T_j T_i \qquad \text{if $ |i-j| > 1$}, \\
y_i y_j & = & y_j y_i , \label{e: yy=yy}\\
y_i T_j & = & T_j y_i \qquad \text{if $ i \neq j,j+1 $}, \\
T_i y_i & = & y_{i+1} T_i^{-1}. 
\eqa 

We consider two different actions of $\Hhat$ in   $ \cplx[z_1^{\pm 1},\dots, z_N^{\pm 1}]$.

The first of these actions is defined by:    
\begin{equation}
y_i = z_i^{-1} ,\quad  \text{and}  \quad  T_i = g_{i,i+1} ,  \label{e: ahecke1}
\end{equation}

where the operators $ g_{i,j} $ are as follows \cite{BG} :
\begin{equation}
g_{i,j} =  \frac{ q^{-1}z_i - q z_j }{ z_i - z_j } ( K_{i,j} - I) + q , \quad 1\leq i\neq j \leq N ,   \label{e: gij}
\end{equation}
and $ K_{i,j} $ is the permutation operator for  variables $z_i$ and $z_j$.

The other action of $\Hhat$ is specified by:
\begin{equation}
y_i = q^{1-N}Y^{(N)}_i , \quad \text{and} \quad  T_i = g_{i,i+1} , \label{e: ahecke2} 
\end{equation}
here $ Y^{(N)}_i $ are difference operators:
\begin{align}
& Y^{(N)}_i   =  K_{i,i+1}g_{i,i+1}^{-1}\dots K_{i,N}g_{i,N}^{-1}p^{D_i}K_{1,i}g_{1,i}\dots K_{i-1,i}g_{i-1,i} , \label{e: Y}\\
\intertext{and}
& p^{D_i}f(z_1,\dots,z_i,\dots,z_N)  =  f(z_1,\dots,pz_i,\dots,z_N) , \quad f \in \cplx[z_1^{\pm 1},\dots, z_N^{\pm 1}]. \nonumber 
\end{align}

Throughout this paper the $q$ is taken to be a complex number which is not a root of unity ($q =1$ is allowed). We call such $q$ and a $p \in \cplx \setminus q^{2 {\Bbb Q}_{\geq 0}}$ generic, and in what follows we consider only generic  $q$ and $p$ unless stated otherwise. 

\subsection{Eigenfunctions of the operators  $Y^{(N)}_i$ } 

To each $\lb := (\l_1,\l_2,\dots,\l_N) $ $\in$ $\zint^N$ corresponds a monomial $ \zz^{\lb} := z_1^{\l_1}z_2^{\l_2}\dots z_N^{\l_N} $ of the total degree $ |\lb | := \l_1+\l_2+\dots + \l_N $. Let $\lb^+$ be the unique partition ( we permit negative parts ) obtained by ordering the elements of $\lb$: $ \lb^+ = (\l^+_1 \geq \l_2^+ \geq \dots \geq \l^+_N)$, $ \l_{\s(i)}^+ = \l_i $ for a suitable permutation $ \s := \{ \s(1),\s(2),\dots,\s(N)  \} $ of the set $\{1,2,\dots,N\}$. We fix such a $ \s$ uniquely by requiring $ \s(i) < \s(j) $ whenever $ i < j $ and $ \l_i = \l_j $.     

There is a partial order relation in $\zint^N$. First, for two partitions $ \lb^+ = (\l^+_1 \geq \l_2^+ \geq \dots \geq \l^+_N)$ and $ \mb^+ = (\mu^+_1 \geq \mu_2^+ \geq \dots \geq \mu^+_N)$ the partial order is defined by: 
\begin{equation}
\lb^+ \succeq \mb^+ \quad \Leftrightarrow \quad \sum_{j=1}^i \l^+_j \geq \sum_{j=1}^i \mu^+_j
\; \quad  (i=1,2,\dots,N),\quad |\lb^+| = |\mb^+|. 
\end{equation}
This order is extended to $\zint^N$ as follows. For $\lb,\mb$ $\in $ $\zint^N$ put  $\lb \succ \mb$ if:
\begin{align}
& \text{either $ \lb^+ \succ \mb^+ $ ,  or $ \lb^+ =  \mb^+ $  and }  \\   
& \text{the last non-zero element in $( \l_1 - \mu_1 , \l_2 - \mu_2,\dots , \l_N - \mu_N )$ is negative. } \nonumber 
\end{align}

It is straightforward to verify that the action of $ \xi_{i,j} :=K_{i,j} g_{i,j} $ and the $Y_i^{(N)}$ in the monomial basis is triangular. More precisely:  
\begin{align}
& \xi_{i,j}\zz^{\lb} = \begin{cases} q^{-1}\zz^{\lb} + \sum_{\lb^+ \succ \mb^+}c(\lb,\mb)\zz^{\mb} & (\l_i < \l_j) , \\
 q \zz^{\lb} &  (\l_i = \l_j) , \\
q\zz^{\lb} + (q-q^{-1}) \zz^{(i,j)\lb}+\sum_{\lb^+ \succ \mb^+}c(\lb,\mb)\zz^{\mb} & (\l_i > \l_j) \end{cases} \tag{i} \label{e: lemma1i} \\ \intertext{where $(i,j)\lb := \lb |_{\l_i \leftrightarrow \l_j}$.}  
& Y_i^{(N)} \zz^{\lb} = p^{\l_i}q^{2 \s(i) - N -1} \zz^{\lb} +  \sum_{\lb \succ \mb }c(\lb,\mb)\zz^{\mb} . \tag{ii} \label{e: lemma1ii}
\end{align}
Let us put $ \zeta_i(\lb) := p^{\l_i}q^{2 \s(i) - N -1} $ for $ \lb \in \zint^N$ and $ i=1,2,\dots,N$. Since for generic $p$ and $q$ the equality $ \zeta_i(\lb) =  \zeta_i(\mb) $ $(i=1,2,\dots,N)$ implies $\lb = \mb $ we immediately come to the conclusion that the $ Y_i^{(N)} $ admit a common  eigenbasis $ \{ \Phi^{\lb}(\zz)\; |\; \lb \in \zint^N \}$:
\begin{align}
Y_i^{(N)}\Phi^{\lb}(\zz) & = \zeta_i(\lb) \Phi^{\lb}(\zz)  , \quad (i=1,2,\dots,N), \label{e: cor1}\\ \intertext{and} 
\Phi^{\lb}(\zz) &  = \zz^{\lb } + \sum_{\lb \succ \mb }c(\lb,\mb)\zz^{\mb}. \label{e: cor2}
\end{align}  

Following \cite{Cherednik} we will refer to  the Laurent polynomials $ \Phi^{\lb}(\zz) $  as Non-symmetric Macdonald Polynomials ( of type $A$ ). 

The action of the finite Hecke algebra generators $g_{i,i+1}$ in the basis  $ \{ \Phi^{\lb}(\zz)\; |\; \lb \in \zint^N \}$ is summarized as follows \cite{Uglov1}: 

\begin{gather}
g_{i,i+1}\Phi^{\lb}(\zz) = A_i(\lb)\Phi^{\lb}(\zz) + B_i(\lb)\Phi^{(i,i+1)\lb}(\zz), \label{e: heckeact1}\\ \text{where $(i,i+1)\lb := \lb |_{\l_i \leftrightarrow \l_{i+1}}$ and:} \nonumber \\ 
A_i(\lb) := \frac{(q-q^{-1})x}{x-1} , \quad B_i(\lb) := \begin{cases} q^{-1}\{x\} & \quad ( \l_i > \l_{i+1} ) ; \\  
 0  & \quad ( \l_i = \l_{i+1} ) ; \\
q^{-1} & \quad ( \l_i < \l_{i+1} ), \end{cases} \label{e: heckeact2}\\
\{x\} := \frac{(x-q^2)(q^2 x-1)}{(x-1)^2} , \quad x := \frac{\zeta_{i+1}(\lb)}{\zeta_i(\lb)}.
\end{gather}
Note that when $i$ is such that $\l_i = \l_{i+1}$ we have $\frac{\zeta_{i+1}(\lb)}{\zeta_i(\lb)} = q^2$ and hence $ g_{i,i+1}\Phi^{\lb}(\zz) = q \Phi^{\lb}(\zz) $.   

\subsection{ Finite $q$-wedge Product }

Let $n\geq 2$ and $N$ be positive integers. We set: $ V:= \cplxn $ with a base $\{v_1,v_2,\dots,v_n\}$ and $V(z):= \cplx [z^{\pm 1}]\otimes V$ with a base $ \{ z^m v_{\ep} \}$, $ m\in \zint $ , $ \ep \in \{1,2,\dots,n\} $. Often it will be convenient to set $ k = \ep - nm $ , $ u_k := z^m v_{\ep} $. Then $\{ u_k \; | \; k \in \zint \} $ is a base in $V(z)$. In what follows we will use both  notations: $ u_k $ and $ z^m v_{\ep} $ switching between them without further alert.      
The $q$-wedge product of spaces $V(z)$ is defined as a suitable quotient of the tensor product $\otimes^N V(z) \cong \cplx[z_1^{\pm 1},z_2^{\pm 1},\dots,z_N^{\pm 1}]\otimes (\otimes^N V) $. To describe this quotient introduce an action of the finite Hecke algebra in $\otimes^N V $:
\begin{align}
T_i & =  S_{i,i+1},  \quad (i=1,2,\dots,N-1) \label{e: hecke}\\
S & =  -q^{-1}\sum_{1\leq \ep\leq n} E^{\ep,\ep}\otimes E^{\ep,\ep}+ (q-q^{-1})\sum_{1\leq \ep < \ep'\leq n} E^{\ep,\ep}\otimes E^{\ep',\ep'}-\sum_{1\leq \ep \neq \ep'\leq n} E^{\ep,\ep'}\otimes E^{\ep',\ep}, \nonumber 
\end{align}
where $ E^{\ep',\ep} \in End(V) $ is specified by $ E^{\ep',\ep}v_{\alpha} = \delta_{\ep,\alpha}v_{\ep'} $ and  $S_{i,i+1}$ signifies $S$ acting in the $i$-th and $i+1$-th factors in $\otimes^N V $ .

\begin{rmk} The Hecke generators $T_i$  that  are  used in \cite{KMS} are related to the generators which we use in this paper as follows: 
\begin{equation}
T_i = q K_{i,i+1} ( g_{i,i+1} + S_{i,i+1}^{-1} ) - I . \label{e: rmk1}
\end{equation}
\end{rmk}
Now define $ \Omega (\subset \otimes^N V(z))$ as:  
\begin{equation} 
\Omega = \sum_{i=1}^{N-1} Ker( g_{i,i+1} + S_{i,i+1}^{-1} ). \label{e: omega}
\end{equation}
In this setting the $q$-wedge product $\wedge^N V(z)$ is defined as the quotient:
\begin{equation}
\wedge^N V(z) = \otimes^N V(z)/\Omega . \label{e: wedgep}
\end{equation}
This definition is equivalent to the definition in \cite{KMS} due to the Remark above. Notice that for $q=1$ the $q$-wedge product is just the usual exterior (wedge) product of the spaces $V(z)$. In what follows we will use the term ``wedge product'' always for the $q$-deformed wedge product (\ref{e: wedgep}).

Let $\Lambda: \otimes^N V(z) \rightarrow \wedge^N V(z) $ be the quotient map specified by (\ref{e: wedgep}). The image of a pure tensor $ u_{k_1}\otimes u_{k_2}\otimes \dots \otimes u_{k_N} $ under this map is called a wedge and is denoted by:
\begin{equation} 
u_{k_1}\wedge u_{k_2}\wedge \dots \wedge u_{k_N} :=\Lambda ( u_{k_1}\otimes u_{k_2}\otimes \dots \otimes u_{k_N}) . \label{e: wedge} 
\end{equation}
In \cite{KMS} it  is proven that a basis in  $\wedge^N V(z)$ is formed by the normally ordered wedges, that is the wedges (\ref{e: wedge}) such that $ k_1 > k_2 > \dots > k_N $. Any wedge can be written as a linear combination of the normally ordered wedges by using the normal ordering rules \cite{KMS}: 
\bqa
u_l\wedge u_m & = & -u_m\wedge u_l, \quad \text{for $ l=m\bmod n $}, \label{e: no1}\\
u_l\wedge u_m & = & -qu_m\wedge u_l + (q^2 - 1)(u_{m-i}\wedge u_{l+i} - qu_{m-n}\wedge u_{l+n} + \nonumber\\
& & \qquad \qquad + q^2 u_{m-n-i}\wedge u_{l+n+i} + \dots ), \label{e: no2}\\ 
& & \text{for $ l < m , m - l =i \bmod n , 0<i<n$}. \nonumber   
\eqa
The sum above continues as long as the wedges in the right-hand side are normally ordered.

\subsection{ Actions of  $U'_q(\sll)$ in the wedge product.} 

 With each of the two $\Hhat$ actions (\ref{e: ahecke1},\ref{e: ahecke2}) we associate an action of the algebra $U'_q(\sll)$ in the tensor product $\otimes^N V(z)$. Let $ E^{\ep,\ep'}_i$ denote the operator which acts trivially in all factors in  $\cplx[z_1^{\pm 1},\dots, z_N^{\pm 1}]\otimes ( \otimes^N \cplxn ) $ except the $i$-th factor in $\otimes^N \cplxn $ where it acts as $E^{\ep,\ep'}$. The action of the generators $\{ E_{\ep}, F_{\ep} ,K_{\ep} \} $, $(\ep \in \{0,1,\dots,n-1\})$ of $U'_q(\sll)$ (our conventions on $U'_q(\sll)$ are summarized in the Appendix) then is defined  as follows:
\begin{align}
& K_{\ep}  =   K^{\ep}_1K^{\ep}_2\dots K^{\ep}_N , \quad K_i^{\ep} = q^{E_i^{\ep,\ep} - E_i^{\ep +1,\ep +1}} , \label{e: Kfin}\\
& E_{\ep}  =  \sum_{i=1}^N y_i^{-\delta_{0,\ep}}E_i^{\ep,\ep + 1}K_{i+1}^{\ep}\dots K_N^{\ep}, \label{e: Efin}\\  
& F_{\ep}  =  \sum_{i=1}^N y_i^{\delta_{0,\ep}}(K_{1}^{\ep})^{-1}\dots (K_{i-1}^{\ep})^{-1}E_i^{\ep +1,\ep} , \quad \ep = 0,\dots , n-1 \;, \label{e: Ffin} \\
& \text{ where in the right-hand side we regard the indices $\ep$,$\ep+1$ modulo $n$. }  \nonumber 
\end{align}
The substitution $ y_i = z_i^{-1} $ in these expressions  gives the $U'_q(\sll)$ action which was considered in \cite{KMS} --  denote this action by $U^{(N)}_1$. The other choice of the affine Hecke algebra generators: $ y_i = q^{1-N}Y_i^{(N)} $ gives another action of $U'_q(\sll)$ in $\otimes^N V(z) $ -- this action the principal object of study in the present paper , we denote it by $U^{(N)}_0$.

The centre of $\widehat{H}_N(q)$ is generated by symmetric polynomials in $y_i^{\pm 1} $ so we consider two abelian algebras: $H^{(N)}_1$ with generators \cite{KMS}: 
\begin{equation}
B^{(N)}_a := z_1^a + z_2^a + \dots + z_N^a , \quad (a \in \{ \pm1,\pm2,\dots,\});
\end{equation}
And $H^{(N)}_0$ with generators \cite{BG}: 
\begin{equation}
h^{(N)}_a := (q^{1-N}Y^{(N)}_1)^a + (q^{1-N}Y^{(N)}_2)^a + \dots + (q^{1-N}Y^{(N)}_N)^a , \quad (a \in \{ \pm1,\pm2,\dots,\}).\label{e: hamiltonians}
\end{equation}
Obviously $H_j^{(N)}$ commutes with $U_j^{(N)}$ for $j=0,1$.

It is straightforward to verify by using the relations of the affine Hecke algebra, and the explicit form of the operator $S$ (\ref{e: hecke}) that the actions $U_j^{(N)}$,$H_j^{(N)}$ for $j=0,1$ preserve the subspace $\Omega$.   

This implies that $U_j^{(N)}$, $H_j^{(N)}$; $(j=0,1)$ are well-defined in the wedge product $\wedge^N V(z)$, and from now on we consider these actions as defined in $\wedge^N V(z)$.

\subsection{Semi-infinite wedge product}

For $M\in \zint$ the space of semi-infinite wedges $F_M$ is defined in \cite{KMS} as the linear span of semi-infinite monomials:
\begin{equation}
 u_{k_1}\wedge u_{k_2}\wedge u_{k_3}\wedge \quad \dots \quad , \label{e: siwedge}
\end{equation}
such that for $ i >> 1 $ the asymptotic condition $ k_i = M-i+1$ holds.
The vacuum semi-infinite monomial  in $F_M$ is specified by $ k_i = M-i+1, $  $ i=1,2,\dots $ and is denoted by $ |M\rangle $:  
\begin{equation}
|M\rangle = u_M\wedge u_{M-1} \wedge u_{M-2} \wedge \quad \dots \quad .  \label{e: vacuum}
\end{equation}
The normal ordering rules (\ref{e: no1}, \ref{e: no2}) imply that the normally ordered semi-infinite monomials -- that is (\ref{e: siwedge}) with $ k_1 > k_2 > k_3 > \dots $ form a basis in $F_M$.

The level-0 action $U_1^{(N)}$ in the limit $N \rightarrow \infty $ was used in \cite{KMS} to define a {\em level-1} action of $U_q'(\sll)$ in the space $F_M$, such that as an $U_q'(\sll)$-module $F_M$ is isomorphic to the Fock space module introduced in \cite{H}. The abelian algebra $H_1^{(N)}$ in the same limit gives rise to an action of the {\em Heisenberg algebra} in $F_M$.        
\\ \mbox{} \\
The two main problems which we address in the present paper are: 1) To define a {\em level-0}  $U_q'(\sll)$-action in $F_M$ starting from the action $U_0^{(N)}$ in $\wedge^N V(z)$. 2) To construct the irreducible decomposition of the Fock space $F_M$ with respect to this action.  

\section{Decomposition of the finite wedge product}

In this section we find the decomposition of the wedge product $\wedge^N V(z)$ with respect to the  $U_q'(\sll)$-action $U_0^{(N)}$. In order to derive this decomposition we construct a suitable base of $\wedge^N V(z)$ by using the Non-symmetric Macdonald Polynomials.

\subsection{ A base in $\wedge^N V(z)$ }

Let $ \ee := (\ep_1,\ep_2,\dots,\ep_N) $ where $ \ep_i \in \{1,2,\dots,n\} $. For a sequence $\ee $ we set
\begin{equation}
{\bold v}_{\ee}:= v_{\ep_1}\otimes v_{\ep_2}\otimes \dots \otimes v_{\ep_N}\quad  ( \in \otimes^N \cplxn ).  
\end{equation}
A sequence $\mm := (m_1,m_2,\dots,m_N)$ from $\zint^N$ is called $n$-strict if it contains no more than $n$ equal elements of any given value. 
Let us define the sets $\MC_N^n$ and $\EC(\mm)$ by 
\begin{align}
 & \MC_N^n := \{\mm = (m_1,m_2,\dots,m_N) \in \zint^N \; | \; m_1\leq m_2 \leq \dots \leq m_N,\; \text{$\mm $  is $n$-strict } \}, \label{e: Mm}\\
 \intertext{ and for $ \mm \in \MC_N^n $ } 
& \EC(\mm) := \{ \ee = (\ep_1,\ep_2,\dots,\ep_N) \; | \; \ep_i > \ep_{i+1} \; \text{for all $i$ s.t. $m_i =m_{i+1}$ } \}.   \label{e: Em}
\end{align}
In these notations the set 
\begin{equation}
\{ w(\mm,\ee):= \Lambda( \zz^{\mm} \otimes {\bold v}_{\ee}) \equiv z^{m_1}v_{\ep_1}\wedge z^{m_2}v_{\ep_2}\wedge\dots\wedge z^{m_N}v_{\ep_N} \quad | \quad \mm \in \MC_N^n , \ee \in \EC(\mm) \}.
\end{equation}
is nothing but the base of the normally ordered wedges in $\wedge^N V(z)$:
\begin{multline}
\{ w(\mm,\ee) \quad | \quad \mm \in \MC_N^n , \ee \in \EC(\mm) \} = \{ u_{k_1}\wedge u_{k_2} \wedge \cdots \wedge u_{k_N} \quad |  \quad k_1 > k_2 > \dots > k_N \}, \\  k_i = \ep_i - n m_i . \nonumber
\end{multline}

 For the purpose of the $U_q'(\sll)$-decomposition we construct another base. The elements of this new base have the same labels  $\mm \in \MC_N^n , \ee \in \EC(\mm) $ as the elements of the base of the normally ordered wedges.

For  $\mm \in \MC_N^n , \ee \in \EC(\mm) $ let us put
\begin{equation}
\phi(\mm,\ee) := \Lambda( \Phi^{\mm}(\zz) \otimes {\bold v}_{\ee}).
\end{equation}
Notice that at $q=1$ we have $ \phi(\mm,\ee) \equiv w(\mm,\ee) $.

\begin{prop}
The set $\{\phi(\mm,\ee) \; | \;  \mm \in \MC_N^n , \ee \in \EC(\mm) \}$ is a base of $\wedge^N V(z)$.  \label{p: base} 
\end{prop}
\begin{pf}
To show that  $\wedge^N V(z)$ = $span_{\cplx}\{\phi(\mm,\ee) \; | \;  \mm \in \MC_N^n , \ee \in \EC(\mm) \}$ we use the formulas (\ref{e: heckeact1}) for the action of the Hecke algebra generators on Non-symmetric Macdonald Polynomials together with the fact that for any $ i\in \{1,2,\dots,N-1\} $ we have 
\begin{equation}
Im(g_{i,i+1} - S_{i,i+1}) \subset Ker(g_{i,i+1} + S^{-1}_{i,i+1}) \Rightarrow  
\Lambda((g_{i,i+1} - S_{i,i+1})f) = 0 \quad \forall f \in \otimes^N V(z). \label{e: intr} 
\end{equation}
First, the eq. (\ref{e: heckeact1}) and (\ref{e: intr}) allow us to write ( at generic $p$ and  $q$):  
\begin{equation}
\wedge^N V(z) = span_{\cplx}\{\Lambda(\Phi^{\lb}(\zz)\otimes {\bold v}) \; | \;  \lb \in \zint^N , \l_1 \leq \l_2 \leq \dots \leq \l_N ; {\bold v}  \in \otimes^N\cplxn \}. \label{e: int0}
\end{equation}
Next, we observe that (\ref{e: heckeact1},\ref{e: intr}) imply:
\begin{equation}
\Lambda ( \Phi^{\lb}(\zz)\otimes ( q + S_{i,i+1}^{-1}) {\bold v} ) = 0 \quad \text{whenever $ \l_i = \l_{i+1}. $} \label{e: int1}
\end{equation}
From the last equation it follows that a vector $\Lambda ( \Phi^{\lb}(\zz)\otimes {\bold v} )$ is equal to zero if the sequence $ \lb $ is not $n$-strict. Now using  (\ref{e: int1}) together with the explicit form of the operators $ S_{i,i+1} $ (\ref{e: hecke}) we derive from  (\ref{e: int0}) that the vectors    $ \phi(\mm,\ee):= \Lambda(\Phi^{\mm}(\zz) \otimes {\bold v}_{\ee} ) $ with $  \mm \in \MC_N^n $ and $ \ee \in \EC(\mm) $ span the wedge product.  

To demonstrate linear independence of these vectors we consider the limit $q=1$ in which limit the $\phi(\mm,\ee)$ coincide with the $w(\mm,\ee)$ -- elements of the base of the normally ordered wedges.
\end{pf}

Since $\Phi^{\mm}(\zz)$ is an eigenvector of the operators $ Y_i^{(N)} $ it is clear, that $ U_0^{(N)} $ and $H_0^{(N)} $ preserve the subspace 
\begin{equation}
 E^{\mm} :=  \oplus_{\ee \in \EC(\mm)} \cplx \phi(\mm,\ee) \label{e: Emm1} 
\end{equation}
for any $\mm \in \MC_N^n $. Moreover it is easy to see, that the $ E^{\mm} $ is an eigenspace of the operators $ h_a^{(N)} $ which generate $H_0^{(N)} $. The eigenvalue $ h_a^{(N)}(\mm) $ of $ h_a^{(N)} $  $ (a=\pm1,\pm2,\dots \;)$ for this eigenspace is
\begin{equation} 
h_a^{(N)}(\mm) = \sum_{i=1}^N ( q^{1-N} \zeta_i(\mm))^a  = \sum_{i=1}^{N} p^{a m_i}q^{2a(1-i)}. 
\label{e: heigenvalueN}\end{equation}
Note also, that the $ E^{\mm} $ is an eigenspace of the degree operator $ z_1\frac{\p}{\p z_1} + z_2\frac{\p}{\p z_2} + \dots + z_N\frac{\p}{\p z_N} $ with the eigenvalue $ |\mm |:= m_1+m_2+\dots +m_N $.

In the rest of this section we will describe the structure of the  $ E^{\mm} $ for a fixed $\mm \in \MC_N^n $ as an $U_q'(\sll)$-module with the action $ U_0^{(N)} $. 

\subsection{ The $U_q'(\sll)$-module $ E^{\mm} $ }

 For $a_1,a_2,\dots,a_r$ $\in \cplx$ let $ \pi_{a_1,\dots,a_r}^{(r)} $ be the evaluation action of $U_q'(\sll)$ defined in $\otimes^N \cplxn$ by
\begin{align}
& \pi_{a_1,\dots,a_r}^{(r)}(K_{\ep})  =   K^{\ep}_1K^{\ep}_2\dots K^{\ep}_r , \quad K_i^{\ep} = q^{E_i^{\ep,\ep} - E_i^{\ep +1,\ep +1}} , \label{e: Kfin2}\\
&\pi_{a_1,\dots,a_r}^{(r)} (E_{\ep})  =  \sum_{i=1}^r a_i^{\delta_{0,\ep}}E_i^{\ep,\ep + 1}K_{i+1}^{\ep}\dots K_r^{\ep}, \label{e: Efin2}\\  
& \pi_{a_1,\dots,a_r}^{(r)}(F_{\ep})  =  \sum_{i=1}^r a_i^{-\delta_{0,\ep}}(K_{1}^{\ep})^{-1}\dots (K_{i-1}^{\ep})^{-1}E_i^{\ep +1,\ep} , \quad \ep = 0,\dots , n-1 \;, \label{e: Ffin2} \\
& \text{ where in the right-hand side we regard the indices $\ep$,$\ep+1$ modulo $n$. }  \nonumber 
\end{align}
In this notation we set for $\mm \in \MC_N^n $:
\begin{equation}
\pi(\mm) := \pi_{a_1,\dots,a_N}^{(N)} \quad \text{wherein we put $a_i = q^{N-1}\zeta_i(\mm)^{-1}$.}  \label{e: pim}
\end{equation}
Since $ \zeta_{i+1}(\mm) = q^2 \zeta_i(\mm) $ whenever $m_i=m_{i+1}$ the action $\pi(\mm)$  has the following property:
\begin{equation}
( q + S_{i,i+1}^{-1})\pi(\mm) =  \pi(\mm) |_{\zeta_{i+1}(\mm) \leftrightarrow\zeta_i(\mm)}( q + S_{i,i+1}^{-1})\quad \text{for all $i$ s.t. $m_i=m_{i+1}$ }.
\end{equation}
As the result the  $\pi(\mm)$ is well-defined in  the space:
\begin{equation}
W^{\mm} := \otimes^N\cplxn / \sum_{\{i | m_i = m_{i+1} \}} Ker( q + S_{i,i+1}^{-1}). \label{e: Wm}
\end{equation}

Let $\iota_{\mm}$: $\otimes^N\cplxn \rightarrow W^{\mm} $ be the quotient map defined by (\ref{e: Wm}). We define a map 
\begin{equation}
\beta :  E^{\mm} = span_{\cplx}\{ \Lambda(\Phi^{\mm}(\zz)\otimes v) \;|\; v \in  \otimes^N\cplxn \} ( = \oplus_{\ee \in \EC(\mm) } \cplx \Lambda(\Phi^{\mm}(\zz)\otimes {\bold v}_{\ee}) ) \; \rightarrow \; W^{\mm}  \label{e: beta}   
\end{equation}
by setting $ \beta( \Lambda(\Phi^{\mm}(\zz)\otimes v) ) $ = $ \iota_{\mm}(v) $.

\begin{prop}
The map $\beta$ \em{(\ref{e: beta})} is an isomorphism of the $U_q'(\sll)$-modules  $( E^{\mm} , \; U_0^{(N)})$  and  $( W^{\mm} , \;\pi(\mm)).$
\end{prop}
\begin{pf}
To show that the $\beta$ is an isomorphism of the linear spaces it is sufficient to observe that the set $ \{ \iota_{\mm}({\bold v}_{\ee}) \; | \; \ee \in \EC(\mm) \} $ is a base in $W^{\mm}$.

The map $\beta$ is an intertwiner of the actions $U_0^{(N)}$ and $ \pi(\mm)$ since for any generator $x$ of $U_0^{(N)}$ we have: 
\begin{equation}
x. \Lambda(\Phi^{\mm}(\zz)\otimes v) = \Lambda(\Phi^{\mm}(\zz)\otimes \pi(\mm)(x).v).
\end{equation}
\end{pf}

Now let us subdivide the sequence $\mm = (m_1\leq m_2 \leq \dots \leq m_N )$ into subsequences which comprise equal elements: 
\begin{equation}
\mm = ( m_1 = \dots  = m_{r_1} < m_{1+r_1} =  \dots = m_{r_2} < \cdots <  m_{1+r_J} = \dots = m_N ). \label{e: splitm}
\end{equation}
Notice that since $\mm$ is $n$-strict we have for $ 1\leq k \leq J+1$ the inequalities $ 1\leq r_k - r_{k-1} \leq n $ (here we put $r_0:=0$ and $r_{J+1}:=N$). As the  $U_q'(\sll)$-module the space $W^{\mm}$ then is represented as the following tensor product: 
\begin{equation}
W^{\mm} = V[p^{-m_{r_1}}q^{2(r_1 - 1)},r_1 ]\otimes V[p^{-m_{r_2}}q^{2(r_2 - 1)},r_2 - r_1]\otimes \dots \otimes V[p^{-m_N}q^{2(N - 1)},N - r_J] \label{e: Wm1}
\end{equation}
where $V[a,j]$ $(1\leq j \leq n)$ is the $U_q'(\sll)$-module with the action $\pi^{(j)}_{a_1,\dots,a_j}$ such that $ a_1:= a ,$ $ a_2:= q^{-2} a ,$ $\dots ,$ $ a_j:= q^{-2(j-1)}a $ (\ref{e: Kfin2} - \ref{e: Ffin2}), and as a  linear space
\begin{equation} 
V[a,j] =  \otimes^j\cplxn / \sum_{i=1}^{j-1} Ker(q+S_{i,i+1}^{-1}).  \label{e:
Vaj}   
\end{equation}

Our next task is to describe the $V[a,j]$ $(1 \leq j \leq n)$. First of all for $ 1\leq j \leq n-1$ as $U_q(\sln)$-module ( $U_q(\sln) \subset U_q'(\sll)$ ) the  $V[a,j]$ is irreducible and isomorphic to the highest weight module $V(\Lambda_j)$ with the fundamental $U_q(\sln)$  highest weight  $\Lambda_j$. When $j=n$; the  $V[a,j]$ is the 1-dimensional trivial representation of $U_q(\sln)$. Thus $V[a,j]$ $(1 \leq j \leq n)$ is an irreducible $U_q'(\sll)$-module, and in order to give the complete specification of $V[a,j]$ for $1 \leq j \leq n-1$ it is sufficient to describe  associated Drinfel'd Polynomials. In the conventions of \cite{CP} which we recall in the Appendix these are provided by the following Lemma:
\begin{lemma}
For $1\leq j \leq n-1$ the $V[a,j]$ is an irreducible $U_q'(\sll)$-module with the Drinfel'd Polynomials
\begin{equation}
P_k(u) = \begin{cases} u - q^{j-2} a^{-1} & \text{\em{for} $k=j$} , \\
                       1 & \text{ \em{for} $1\leq k \leq n-1$, $k\neq j$}. \end{cases}
\end{equation}  \label{l: Drinfeld}
\end{lemma}

Using the results of  \cite{AK} we can claim, that for generic $q$ and $p$ the representation $W^{\mm}$ (\ref{e: Wm1}) is irreducible, therefore the Drinfel'd Polynomials of $W^{\mm}$ are just products of the Drinfel'd Polynomials associated with the factors  $ 
V[p^{-m_{r_k}}q^{2(r_k - r_{k-1})},r_k-r_{k-1}]$  (Cf.\cite{CP}). This leads to the main Proposition of this section:
\begin{prop}
For generic $p$ and $q$ the $E^{\mm} \cong W^{\mm}$ is irreducible and with notations of \em{(\ref{e: splitm})} the Drinfel'd Polynomials of $E^{\mm} \cong W^{\mm}$ are
\begin{equation}
P_i(u) =  \prod_{\{1\leq k \leq J+1\; | \; r_k - r_{k-1} = i \}} (u - p^{m_{r_k}}q^{-r_k-r_{k-1}}) \qquad ( i\in \{ 1,2,\dots,n-1 \}). 
\end{equation}  \label{p: finwedgedecomp}
\end{prop}
The proofs of this Proposition and the Lemma \ref{l: Drinfeld} are discussed in the Appendix.

Finally we note that 
\begin{equation}
\wedge^N V(z) = \bigoplus_{\mm \in \MC_N^n} E^{\mm}. 
\end{equation}
This is the desired decomposition of the finite wedge product with respect to the $U_q'(\sll)$ action $U_0$.

\section{ A level-0 action of $U_q'(\sll)$ in the Fock space } 

In this section we will define a level-0 action of the $U_q'(\sll)$ in the space of semi-infinite wedges $F_M$ ( or equivalently in the level-1 Fock space module of the same algebra -- $U_q'(\sll)$ --  see subsection 2.5 )  starting from the action $U_0^{(N)}$ which was defined in sec. 2 in the finite wedge product. The level-0 action in the space $F_M$ is constructed by taking a suitable projective limit of  $U_0^{(N)}$ and can be thought of as an appropriate, well-defined limit of   $U_0^{(N)}$ when the number of particles $N$ goes to infinity.

In this section we allow the parameter $p$ to  be arbitrary complex number, whrereas the $q$ is still required to be generic.

 Let $w = u_{k_1}\wedge u_{k_2} \wedge \dots \wedge u_{k_N}$ $\equiv$ $z^{m_1}v_{\ep_1}\wedge z^{m_2}v_{\ep_2} \wedge \dots \wedge z^{m_N}v_{\ep_N}$, $ w\in \wedge^N V(z)$ be a normally ordered wedge: $ k_1 > k_2 > \dots > k_N $. Often it will be convenient to label this wedge by the two sequences $\mm \in \MC_N^n$ and $\ee \in \EC(\mm)$ : 
\bqa 
\mm & = & (m_1 \leq m_2 \leq \dots \leq m_N ), \quad m_i \in \zint , \label{e: mseq}\\
\ee & = & (\ep_1,\ep_2,\dots,\ep_N), \quad  1\leq \ep_i\leq n , \label{e: eseq}\\ 
& & \text{such that $ k_i = \ep_i - nm_i $}, 
\eqa
and write: $ w := w(\mm,\ee) $.
\mbox{} From now on  we  will use the notation $w(\mm,\ee)$ {\em exclusively} for normally ordered wedges.

Recall that by the definitions (\ref{e: Mm}, \ref{e: Em}) of the $\MC_N^n$ and $\EC(\mm)$ the  $\mm$ in (\ref{e: mseq}) is always $n$-strict sequence and that  $ \ep_i > \ep_{i+1} $ whenever $m_i = m_{i+1}$ .

Throughout this section we fix an integer $M$ and  $ 0 \leq s \leq n-1 $ such that  $M = s\bmod n.$  Let $ r \in \{0,1,2,\dots \; \}$; and let $ w(\mm^0,\ee^0) $ $\in $ $ \wedge^{s+nr}V(z)$ be defined as:
\begin{equation}
w(\mm^0,\ee^0):= u_M\wedge u_{M-1} \wedge \dots \wedge u_{M-(s+nr)+1} . \label{e: vacuumN}
\end{equation}
Here the sequences $\mm^0$ and $\ee^0$ are as follows:
\bqa
\mm^0 = (m_1^0,m_2^0,\dots,m_{s+nr}^0) & := & (\underbrace{m^0,\dots,m^0}_s,\underbrace{m^0+1,\dots,m^0+1}_n ,\nonumber \\ & & \underbrace{m^0+2,\dots,m^0+2}_n, \dots ,\underbrace{m^0+r,\dots,m^0+r}_n )  \label{e: m0N} \\
\ee^0 = (\ep_1^0,\ep_2^0,\dots,\ep_{s+nr}^0) & := & (\underbrace{s,s-1,\dots,1}_s,\underbrace{n,n-1,\dots,1}_n ,\nonumber \\ & & \underbrace{n,n-1,\dots,1}_n, \dots ,\underbrace{n,n-1,\dots,1}_n )  \label{e: e0N} , \\
& & \text{where}\quad m^0 := \frac{s-M}{n} \qquad ( 0\leq s\leq n-1).
\eqa
We will call these two sequences {\em vacuum} sequences, and the  $ w(\mm^0,\ee^0) $ -- {\em vacuum vector} of  $ \wedge^{s+nr}V(z)$.

Define $ V_M^{s+nr} \subset \wedge^{s+nr}V(z) $ as:
\begin{equation}
V_M^{s+nr}  = \bigoplus\begin{Sb} \mm \in \MC_N^n, \ee \in \EC(\mm) \\ m_{s+nr} \leq m_{s+nr}^0 \end{Sb} \cplx w(\mm,\ee) . \label{e: VM}
\end{equation}
Notice that the condition $m_{s+nr} \leq m_{s+nr}^0$ in this definition is equivalent to the condition:
\begin{equation}
m_i \leq m_i^0 \quad \text{ for all $ i=1,2,\dots,s+nr$} \label{e: mlessm0}
\end{equation}
because the sequence $\mm$ is $n$-strict and non-decreasing.

\begin{prop} 
The $U_q'(\sll)$ action $U_0^{(s+nr)}$ and the operators $h_l^{(s+nr)}$ $(l\in\zint_{\neq 0})$  preserve the subspace $ V_M^{s+nr} $. \label{p: preserve}
\end{prop}

To prove this and other Propositions we need two Lemmas. The first of these Lemmas concerns properties of the operators (\ref{e: Y}) and operators:
\begin{equation}
\xi_{i,j} = K_{i,j} g_{i,j},
\end{equation}
where $K_{i,j}$ and $g_{i,j}$ are  defined in (\ref{e: gij}).

\begin{lemma}
Let $\zz^{\nn}$ $\equiv$ $z_1^{n_1}z_2^{n_2}\dots z_N^{n_N}$ be a monomial in $\cplx [z_1^{\pm},\dots,z_N^{\pm}]$, \\ and let $a =\max\{n_1,\dots,n_N\}$. 

Then  
\bqa
\xi_{i,j}^{\pm 1}\zz^{\nn} & = &  \sum_{\nn'} c_{\pm}(\nn,\nn')\zz^{\nn'}, \label{e: st1} \\
 (Y_i^{(N)})^k\zz^{\nn} & =  & \sum_{\nn'} c_k(\nn,\nn')\zz^{\nn'}, \quad k = \pm1,\pm2,\dots \quad . \label{e: st2}
\eqa
where $c_{\pm}(\nn,\nn'),c_k(\nn,\nn')$ are coefficients, and the summation ranges over $\nn'$ such that:
\bqa 
 n_1',n_2',\dots,n_N' & \leq & a , \label{e: lemma1a}\\
\#\{n_i' | n_i' = a\} & \leq & \#\{n_i | n_i = a\},\label{e: lemma1b} \\
n_1'+ n_2'+\dots + n_N' & = & n_1+ n_2+\dots + n_N .\label{e: lemma1c}
\eqa \label{l: lemma1}
\end{lemma}
\begin{pf}
To prove the statement about the summation range (\ref{e: lemma1a} - \ref{e: lemma1c}) in (\ref{e: st1}) we use the explicit formulas for the action of $\xi_{i,j}$ and $\xi_{i,j}^{-1}$ on monomials:
\bqa 
 \text{for} & n_i < n_j & : \label{e: fa}  \\ 
\xi_{i,j}^{\pm 1} z_i^{n_i} z_j^{n_j} & = & q^{\mp 1}z_i^{n_i} z_j^{n_j} \mp (q-q^{-1})\sum_{k=1}^{n_j-n_i-1}z_i^{n_i+k} z_j^{n_j-k} , \nonumber \\
 \text{for} & n_i = n_j & : \label{e: fb}  \\ 
\xi_{i,j}^{\pm 1} z_i^{n_i} z_j^{n_j} & = & q^{\pm 1}z_i^{n_i} z_j^{n_j}  , \nonumber \\
 \text{for} & n_i > n_j & :  \label{e: fc} \\ 
\xi_{i,j}^{\pm 1} z_i^{n_i} z_j^{n_j} & = & q^{\pm 1}z_i^{n_i} z_j^{n_j} \pm (q - q^{-1})z_i^{n_j} z_j^{n_i} \pm (q-q^{-1})\sum_{k=1}^{n_i-n_j-1}z_i^{n_i-k} z_j^{n_j+k} \nonumber . 
\eqa
The statements (\ref{e: st1}, \ref{e: lemma1a} - \ref{e: lemma1c} ) immediately follow from these formulas. 
The statements (\ref{e: st2}, \ref{e: lemma1a} - \ref{e: lemma1c} ) follow from (\ref{e: st1}, \ref{e: lemma1a} - \ref{e: lemma1c}) and (\ref{e: Y}).
\end{pf}

We will also need the following Lemma which shows triangularity of the normal ordering. 

\begin{lemma}
Let $v \in \otimes^N \cplxn $, and $\Lambda$ be the quotient map defined by {\em (\ref{e: wedgep})}.

Then in the  notations of Lemma \ref{l: lemma1} the following holds:
\begin{equation}
\Lambda (\zz^{\nn}\otimes v) = \sum_{\nn',\ee} c(\nn,v;\nn',\ee) w(\nn', \ee) ,  
\end{equation}
where $c(\nn,v;\nn',\ee)$ is a coefficient and the summation ranges over $\nn'$ such that $(\nn')^+ \preceq (\nn)^+ $ and consequently:
\bqa 
 n_1',n_2',\dots,n_N' & \leq & a (:= \max\{n_1,\dots,n_N\}) ,\label{e: lemma2a}  \\
\#\{n_i' | n_i' = a\} & \leq & \#\{n_i | n_i = a\},\label{e: lemma2b} \\
n_1'+ n_2'+\dots + n_N' & = & n_1+ n_2+\dots + n_N \label{e: lemma2c}
\eqa
\end{lemma} \label{l: lemma2}
\begin{pf}
Use the normal ordering rules (\ref{e: no1}, \ref{e: no2}).
\end{pf}
\\ \mbox{} \\
\begin{ppf}{\ref{p: preserve}}
Let $w(\mm,\ee)$ $\in$ $V_M^{s+nr}$ -- that is $m_{s+nr} \leq m_{s+nr}^0$.
We prove the Proposition by considering the action of the generators of $U_0^{(s+nr)}$ on $w(\mm,\ee)$. 

Action of any of the generators (\ref{e: Kfin} - \ref{e: Ffin}) of the $U_q(\sln)$ subalgebra of $U_0^{(s+nr)}$ on this wedge results in a linear combination of wedges $w(\mm,\ee')$ with the same sequence $\mm$ as in  $w(\mm,\ee)$. Therefore the action of the $U_q(\sln)$ subalgebra preserves $V_M^{s+nr}$. By the same token $K_0$ (\ref{e: Kfin}) preserves $V_M^{s+nr}$ as well.

Consider now the vector:
\begin{equation}
F_0.w(\mm,\ee) =  \sum_{i=1}^N \Lambda( q^{1-N}Y_i^{(N)}.\zz^{\mm}\otimes (K_{1}^{0})^{-1}\dots (K_{i-1}^{0})^{-1}E_i^{1,n}.{\bold v}_{\ee}),
\end{equation}
where we put $ N := s+nr$ and $ {\bold v}_{\ee}:= v_{\ep_1}\otimes v_{\ep_2} \otimes \dots \otimes v_{\ep_N}$ $\in$ $ \otimes^N \cplxn$. In each of the summands apply first  Lemma \ref{l: lemma1} to express $Y_i^{(N)}.\zz^{\mm}$ as a linear combination of monomials, then apply Lemma \ref{l: lemma2} to express the result as a linear combination of the normally ordered wedges:
\begin{equation}
F_0.w(\mm,\ee) =  \sum_{\mm',\ee'} c(\mm,\ee;\mm',\ee')w(\mm',\ee').
\end{equation}
Due to (\ref{e: lemma1a}) in Lemma \ref{l: lemma1} and (\ref{e: lemma2a}) in Lemma \ref{l: lemma2} in the last formula we have: $ m_{s+nr}'\leq m_{s+nr} \leq m_{s+nr}^0 $ and thus $F_0.w(\mm,\ee)$ $\in$ $V_M^{s+nr}$ by the definition (\ref{e: VM}) of the $V_M^{s+nr}$ .

For the generator $E_0$ and the operators $h_l^{(s+nr)}$ the proof is done by the same arguments as for $F_0$.
\end{ppf}
\\ \mbox{} \\
Introduce the degree $|w(\mm,\ee)|$ of a wedge $w(\mm,\ee) \in V_M^{s+nr}$ by: 
\begin{equation}
|w(\mm,\ee)|  =  \sum_{i=1}^{s+nr} m_i^0 - m_i , \label{e: degree}
\end{equation}
where the vacuum sequence $\mm^0$ is defined in (\ref{e: m0N}).  
We have:
\bqa
V_M^{s+nr} & = & \bigoplus_{k \geq 0} V_M^{s+nr,k}, \\ V_M^{s+nr,k} & = & span_{\cplx}\{w(\mm,\ee) \in V_M^{s+nr}\;|\; |w(\mm,\ee)| = k \}. 
\eqa
The statements (\ref{e: lemma1c}) and (\ref{e: lemma2c}) in Lemmas \ref{l: lemma1} and \ref{l: lemma2} imply that the action $U^{(s+nr)}_0$ and the commuting Hamiltonians preserve the degree:
\begin{equation}
 U_0^{(s+nr)} , \; h_l^{(s+nr)} : V_M^{s+nr,k} \rightarrow V_M^{s+nr,k} , \quad l\in \zint_{\neq 0}, \; k =0,1,\dots \quad . \label{e: preserve}
\end{equation}
Notice, that when $ M = 0\bmod n $ (i.e  $s=0$)  the space  $V_M^{s+nr,0}$ is one-dimensional with the basis $w(\mm^0,\ee^0)$.
\\
\mbox{}
\\
Our main technical tool in defining the level-0 action in the space of semi-infinite wedges is the projection map $\ro$ : 
\begin{equation}
\ro : V_M^{s+nr+n} \rightarrow V_M^{s+nr} , \quad r =0,1,\dots \quad , 
\end{equation}
which we define by specifying its action on the normally ordered wedges as follows: 

Let $w(\mm,\ee)$ $\in$ $V_M^{s+nr+n}$ and let  
\bqa
\mm & = & (m_1,m_2,\dots,m_{s+nr},m_{s+nr+1},\dots,m_{s+nr+n} ), \\
\ee & = & (\ep_1,\ep_2,\dots,\ep_{s+nr},\ep_{s+nr+1},\dots,\ep_{s+nr+n} )
\eqa
be the $\mm$ and $\ee$ sequences labeling the wedge $w(\mm,\ee)$ . 

Remove from $\mm$ and $\ee$ the  last $n$ elements, and denote the obtained sequences by $\mm'$ and $\ee'$: 
\begin{align}
\mm' & =  (m_1,m_2,\dots,m_{s+nr}), \label{e: mprime}\\ 
\ee' & =  (\ep_1,\ep_2,\dots,\ep_{s+nr}), \label{e: eprime}\\
\intertext{so that} 
    &  w(\mm',\ee') \in V_M^{s+nr}.  \nonumber 
\end{align}
The action of $\ro$ is then defined by:
\bqa
\ro.w(\mm,\ee) & = & \begin{cases} w(\mm',\ee') & \text{if}  \;  m_{s+nr+1} = m_{s+nr+2} = \dots = m_{s+nr+n} = \\ &  m_{s+nr+1}^0 = m_{s+nr+2}^0 = \dots = m_{s+nr+n}^0 , \label{e: ro} \\
    0 & \text{otherwise} . \end{cases}
\eqa

\begin{prop} The following holds: 
\begin{align}
& \text{The map $\ro$ preserves the degree :} \qquad  \ro : V_M^{s+nr+n,k}   \rightarrow   V_M^{s+nr,k} , \tag{i} \\
\intertext{ and for all $ k=0,1,\dots $ the   map $\rk := \ro |_{V_M^{s+nr+n,k}}$  is surjective. } 
& \qquad \qquad   \text{For $ k \leq r $  the map $ \rk $ is bijective. } \tag{ii}
\end{align} \label{p: prop2}
\end{prop}
\begin{pf} 
The part (i) follows immediately from the definition of the degree (\ref{e: degree}) and from (\ref{e: ro}).  

To prove the part (ii) let us demonstrate,  that if $w(\mm,\ee) \in  V_M^{s+nr+n} $ is such that $ |w(\mm,\ee)| \leq r $, then:  
\begin{multline}
m_{s+nr+1} = m_{s+nr+2} = \dots = m_{s+nr+n} = \\ =  m_{s+nr+1}^0 = m_{s+nr+2}^0 = \dots = m_{s+nr+n}^0 . \label{e: equality}
\end{multline}

Suppose the last equality does not hold. Then we necessarily have:   
\begin{equation}
m_{s+nr+1} = m_{s+nr+1}^0 - t_0, \quad \text{where $t_0 \geq 1$}. \label{e: je} 
\end{equation}
Since the sequence $\mm$ is $n$-strict and non-decreasing, we  also have:      
\begin{align*}
 & m_{s+nr-n+1} < m_{s+nr + 1},\quad m_{s+nr-n+1}^0 = m_{s+nr + 1}^0 - 1  \Rightarrow  \\
            &  m_{s+nr-n+1} = m_{s+nr-n+1}^0 - t_1, \quad  t_1 \geq t_0 , \\
\intertext{and in general:} 
& m_{s+nr-nl+1} = m_{s+nr-nl+1}^0 - t_l , \quad t_{l} \geq t_{l-1} , \quad l=1,2,\dots,r. 
\end{align*}
Summing up the last equations for $ l=1,2,\dots,r $ and (\ref{e: je}) we find:
\begin{equation*}
|w(\mm,\ee)| = \sum_{i=1}^{s+nr+n} m^0_i - m_i \geq \sum_{l=0}^{r} t_l \geq t_0(r+1) \geq  r+1 .
\end{equation*}
This contradicts $ |w(\mm,\ee)| \leq r $, and therefore (\ref{e: equality}) holds. Taking (\ref{e: ro}) into account we find that  
\begin{equation} Ker(\ro |_{V_M^{s+nr+n, k}}) = 0 \end{equation} when $ k \leq r $. 
\end{pf}
\\
\mbox{}
\\
An important property of the map $\ro$ (\ref{e: ro}) is that this map intertwines the $U_q'(\sll)$-action $ U_0^{(s+nr+n)} $ defined in $ V_M^{s+nr+n} $ $\subset $ $\wedge^{s+nr+n}V(z)$ with the $U_q'(\sll)$-action $ U_0^{(s+nr)} $ defined in $ V_M^{s+nr} $ $\subset $ $\wedge^{s+nr}V(z)$. The map $\ro$ also intertwines the actions of the operators $ h_l^{(N)} - h_l^{(N)}(\mm^0)I$ for $N=s+nr+n$ and $N=s+nr$ -- in this case one needs to  redefine $h_l^{(N)}$ (\ref{e: hamiltonians}) by subtracting the eigenvalue associated with the vacuum sequence (\ref{e: heigenvalueN}) . We summarize this as the  Proposition: 
\begin{prop} For $r=0,1,2,\dots$ the following intertwining relations hold :
\begin{align}
&   \qquad \ro U_0^{(s+nr+n)}  =   U_0^{(s+nr)}\ro , \tag{i} \\
&  \qquad  \ro g_l^{(s+nr+n)}  =   g_l^{(s+nr)} \ro ,\quad l\in \zint_{\neq 0} \tag{ii} \\
& \text{where $ g_l^{(N)} = h_l^{(N)} - h_l^{(N)}(\mm^0)I $,} \nonumber \\
& \text{and $ h_l^{(N)}(\mm^0) = \sum_{j=1}^N p^{lm_{j}^0} q^{2 l (1-j)} $.} \nonumber 
\end{align} \label{p: prop3}
\end{prop}
To prove this Proposition we will need one more Lemma on properties of the operators (\ref{e: Y}): 
\begin{lemma}
Let $ \mm $ = $( m_1,m_2,\dots,m_N) \in \zint^N $ be a sequence such, that:
\begin{equation} 
 m_1,m_2,\dots,m_{N-k} < m_{N-k+1} = m_{N-k+2}=\dots = m_N \; \equiv \; m ; \quad 1 \leq k \leq N.\end{equation}
Then for $ a = \pm1,\pm2,\dots, $ the following relations hold:
\begin{align}
\text{\em{for} $0\leq l \leq k-1$} \qquad  (Y_{N-l}^{(N)})^a \zz^{\mm} & =  p^{am} q^{a(2k-2l-N-1)} \zz^{\mm} + [\dots] \; , \label{e: ia}\\ 
\text{\em{for} $1\leq i \leq N-k$} \qquad  (Y_i^{(N)})^a \zz^{\mm} & =  q^{ak} (Y_i^{(N-k)})^a \zz^{\mm} + [\dots] \; , \label{e: ib}\\ 
\intertext{where $[\dots]$ signifies a linear combination of monomials $\zz^{\nn}$ $\equiv$ $ z_1^{n_1}z_2^{n_2}\dots z_N^{n_N}$ such that:}
n_1,n_2,\dots,n_N & \leq  m ,  \\
\intertext{and}
\#\{n_i | n_i = m\} & <  k.
\end{align} \label{l: lemma3}
\end{lemma}

\begin{pf} Consider first the expression:
\begin{multline}
Y_{N-l}^{(N)} \zz^{\mm} = \xi_{N-l,N-l+1}^{-1}\dots \xi_{N-l,N}^{-1}p^{D_{N-l}}\xi_{1,N-l}\dots \xi_{N-k,N-l}\cdot \\ \cdot\xi_{N-k+1,N-l}\dots \xi_{N-l-1,N-l}\zz^{\mm}, \quad 0\leq l\leq k-1.
\end{multline} 
The eq. (\ref{e: fb}) gives:
\begin{equation}
Y_{N-l}^{(N)} \zz^{\mm} = q^{k-l-1}\xi_{N-l,N-l+1}^{-1}\dots \xi_{N-l,N}^{-1}p^{D_{N-l}}\xi_{1,N-l}\dots \xi_{N-k,N-l}\zz^{\mm}, \quad 0\leq l\leq k-1.
\end{equation} 
In the last expression apply $\xi_{N-k,N-l}$ to the monomial $\zz^{\mm}$ using the formula (\ref{e: fa}): 
\begin{multline}
Y_{N-l}^{(N)} \zz^{\mm} =  q^{k-l-1}\xi_{N-l,N-l+1}^{-1}\dots \xi_{N-l,N}^{-1}p^{D_{N-l}}\xi_{1,N-l}\dots \xi_{N-k-1,N-l}( q^{-1} \zz^{\mm} + [\dots]) , \\  \quad 0\leq l\leq k-1,
\end{multline} 
where the meaning of $[\dots] $ is the same as in the statement of the Lemma.

Now apply $\xi_{i,N-l}$ repeatedly for $ i = N-k-1,N-k-2,\dots,1$ using at each step (\ref{e: lemma1a}, \ref{e: lemma1b}) in Lemma \ref{l: lemma1} to show that $\xi_{i,N-l}([\dots]) = ([\dots])$, and using also  (\ref{e: fa}). This gives: 
\begin{equation}
Y_{N-l}^{(N)} \zz^{\mm} =  q^{k-l-1}\xi_{N-l,N-l+1}^{-1}\dots \xi_{N-l,N}^{-1}( p^m q^{k-N} \zz^{\mm} + [\dots]) , \quad 0\leq l\leq k-1.
\end{equation} 
Lemma \ref{l: lemma1} and (\ref{e: fb}) applied in the last formula yield (\ref{e: ia}) for $ a=1 $. 

Consider now the expression:
\begin{equation}
Y_i^{(N)} \zz^{\mm} = \xi_{i,i+1}^{-1}\dots \xi_{i,N-k}^{-1}\xi_{i,N-k+1}^{-1}\dots \xi_{i,N}^{-1}p^{D_i}\xi_{1,i}\dots \xi_{i-1,i}\zz^{\mm}, \quad 1\leq i\leq N-k.
\end{equation}
Write:
\begin{equation}
p^{D_i}\xi_{1,i}\dots \xi_{i-1,i}\zz^{\mm} = (p^{D_i}\xi_{1,i}\dots \xi_{i-1,i}z_1^{m_1}z_2^{m_2}\dots z_{N-k}^{m_{N-k}}) z_{N-k+1}^m\dots z_{N}^m,
\end{equation}
and observe that due to Lemma \ref{l: lemma1} the expression:
\begin{equation}
p^{D_i}\xi_{1,i}\dots \xi_{i-1,i}z_1^{m_1}z_2^{m_2}\dots z_{N-k}^{m_{N-k}}
\end{equation}
is a linear combination of monomials $ z_1^{n_1}z_2^{n_2}\dots z_{N-k}^{n_{N-k}}$ such that $ n_1,n_2,\dots,n_{N-k} < m $,
and therefore the formula (\ref{e: fa}) implies that:
\begin{align}
& \xi_{i,N}^{-1} p^{D_i}\xi_{1,i}\dots \xi_{i-1,i}\zz^{\mm} = 
 q\, p^{D_i}\xi_{1,i}\dots \xi_{i-1,i}\zz^{\mm} + [\dots] \;. \\
\intertext{Continuing to apply Lemma \ref{l: lemma1} together with (\ref{e: fa}) we get :}
&\xi_{i,N-k+1}^{-1}\dots \xi_{i,N}^{-1} p^{D_i}\xi_{1,i}\dots \xi_{i-1,i}\zz^{\mm} = q^k\, p^{D_i}\xi_{1,i}\dots \xi_{i-1,i}\zz^{\mm} + [\dots]\; .
\end{align}
Finally we act by $\xi_{i,i+1}^{-1}\dots \xi_{i,N-k}^{-1}$ on the last expression and using  Lemma \ref{l: lemma1} to show that $\xi_{i,i+1}^{-1} \dots \xi_{i,N-k}^{-1}([\dots])$ = $([\dots])$ ,  arrive at:
\begin{equation}
Y_i^{(N)} \zz^{\mm} = q^k \,\xi_{i,i+1}^{-1}\dots \xi_{i,N-k}^{-1}p^{D_i}\xi_{1,i}\dots \xi_{i-1,i}\zz^{\mm} + [\dots] \; , \quad 1\leq i\leq N-k 
\end{equation}
which is the statement (\ref{e: ib}) of the Lemma at $ a = 1$. 

For $a=-1$ the proof is completely analogous, and (\ref{e: ia},\ref{e: ib}) for $a = \pm 2,\pm3,\dots $ follow from the case $a=\pm1$ and  Lemma \ref{l: lemma1}, eqs. (\ref{e: lemma1a}, \ref{e: lemma1b}). 
\end{pf}
\\ \mbox{} \\
\begin{ppf}{\ref{p: prop3}}
To prove the part (i) we consider the action of the generators of $U_0^{(s+nr+n)}$ on a wedge $ w(\mm,\ee)$ $\in$ $V_M^{(s+nr+n)}$. 

First let $w(\mm,\ee)$ be such that:
\begin{equation}
 \text{$ m_i < m_i^0 $ for at least one $ s+nr < i \leq s+nr + n $.} \label{e: co}
\end{equation}
This condition is equivalent to $ w(\mm,\ee) \in Ker\ro $. The Lemmas \ref{l: lemma1} and \ref{l: lemma2}  imply that acting with any of the generators (\ref{e: Kfin} -- \ref{e: Ffin}) (where $ y_i = q^{1-(s+nr+n)}Y_i^{(s+nr+n)}$) on such $ w(\mm,\ee)$ produces a linear combination of wedges that have the property (\ref{e: co}) as well, and therefore vanish when acted on by $ \ro $. We formulate this as:
\begin{equation}  
U_0^{(s+nr+n)}: Ker\ro \rightarrow Ker\ro.
\end{equation}
Thus for $w(\mm,\ee)$ that satisfy (\ref{e: co}) and any generator $x$ of $U_q'(\sll)$ one has:
\begin{equation}
x^{(s+nr)}\ro w(\mm,\ee) = \ro x^{(s+nr+n)} w(\mm,\ee) = 0.
\end{equation}

Now let $w(\mm,\ee)$ be such that:
\begin{equation}
 \text{$ m_i = m_i^0 $ for all $ s+nr < i \leq s+nr + n $ .} \label{e: coin}
\end{equation}
With the same $\mm'$ and $\ee'$ as in (\ref{e: mprime}, \ref{e: eprime}) one can write:
\begin{align}
& w(\mm,\ee)   =  w(\mm',\ee')\wedge(u_{M-s-nr}\wedge u_{M-s-nr-1}\wedge \dots \wedge u_{M-s-nr-n+1}), \label{e: xx} \\
& \ro w(\mm,\ee)  =  w(\mm',\ee').\label{e: xxo} 
\end{align}

Apply $F_0^{(s+nr+n)}$ to $w(\mm,\ee)$:
\begin{equation}
F_0^{(s+nr+n)}w(\mm,\ee) = \Lambda(\sum_{j=1}^{s+nr+n} q^{1-(s+nr+n)}Y_j^{(s+nr+n)}\zz^{\mm}\otimes
(K_1^0)^{-1}\dots(K_{j-1}^0)^{-1}E^{1,n}_j {\bold v}_{\ee}).
\end{equation}
Lemma \ref{l: lemma3} for $ a=1$ ,$N=s+nr+n, k=n, m=m_{s+nr+n}^0$ and Lemmas \ref{l: lemma1}, \ref{l: lemma2}  enable us to transform the right-hand side of the last equation and arrive at:
\begin{multline}
F_0^{(s+nr+n)}w(\mm,\ee) = (F_0^{(s+nr)}w(\mm',\ee'))\wedge(u_{M-s-nr}\wedge u_{M-s-nr-1}\wedge \dots \wedge u_{M-s-nr-n+1}) + \\ + p^{m_{s+nr+r}^0}q^{2(1-nr-n-s)}\,((K_0^{(s+nr)})^{-1}w(\mm',\ee'))\wedge \cdot \\ \cdot \wedge (u_{M-s-nr-n+1}\wedge u_{M-s-nr-1}\wedge u_{M-s-nr-2}\wedge  \dots \wedge u_{M-s-nr-n+1}) +  \tilde{w}, \label{e: lst}
\end{multline}
where $\tilde{w} \in Ker\ro$.
The normal ordering rules (\ref{e: no1},\ref{e: no2}) imply that the second summand in the right-hand side of the last equation vanishes ( cf. Lemma 2.2 in \cite{KMS}).
Finally (\ref{e: lst}) and (\ref{e: xx},\ref{e: xxo}) give:
\begin{equation}
F_0^{(s+nr)}\ro w(\mm,\ee) = \ro  F_0^{(s+nr+n)}w(\mm,\ee).
\end{equation}
Thus (i) is proven for the generator $F_0$. The rest of the $U_q'(\sll)$-generators and the statement (ii) of the Proposition are handled in the same way. 
\end{ppf}
\\ \mbox{} \\
At fixed $M=s\bmod n$ and fixed degree $k$ form the projective limit of the spaces $V_M^{s+nr,k}$ with respect to the map $\rk := \ro |_{V_M^{s+nr+n,k}}$: 
\begin{equation}
V_M^k = \lim\begin{Sb} \longleftarrow \\ r \end{Sb} V_M^{s+nr,k}. \label{e: inverselim} 
\end{equation}
A vector in $V_M^k$ is a semi-infinite sequence $\{ f_r \}_{r\geq 0}$ ; $ f_r $ $\in$ $V_M^{s+nr,k}$ such that:
\begin{equation}
\rk f_{r+1} = f_r, \quad r = 0,1,2,\dots \;.
\end{equation}
Since the map $\rk$ is bijective when $r\geq k$ we have the isomorphism of linear spaces: 
\begin{equation}
 V_M^k \cong  V_M^{s+nr,k},\quad r \geq k.
\end{equation}
Notice that for $f = \{ f_r \}_{r\geq 0}$ $\in$ $ V_M^{s+nr,k}$, $ r \geq k$:
\begin{equation}
f_{r+1} = f_r \wedge u_{M-s-nr}\wedge u_{M-s-nr-1}\wedge \dots \wedge u_{M-s-nr-n+1}, \label{e: tail}
\end{equation}
as implied by Proposition \ref{p: prop2} (ii) and the definition (\ref{e: ro}) of the map $\ro$. 

Now we use Propositions \ref{p: prop2} (i) and \ref{p: prop3} to define in  the space $V_M^k $ an $U_q'(\sll)$-action $U_0^{(\infty)}$, and an action of the commutaive family $\{g_l^{(\infty)}\}_{l \in \zint_{\neq 0}}.$  For $ f = \{ f_r \}_{r\geq 0}$ $\in$ $V_M^k $ set:
\bqa
U_0^{(\infty)}.\{ f_r \} & := & \{ U_0^{(s+nr)}.f_r \} , \\
g_l^{(\infty)}.\{ f_r \} & := & \{g_l^{(s+nr)}.f_r \}.
\eqa
Clearly we still have:
\begin{equation}
U_0^{(\infty)}g_l^{(\infty)} = g_l^{(\infty)}U_0^{(\infty)},\quad l \in \zint_{\neq 0}. 
\end{equation}
\\ \mbox{} \\
The degree of a semi-infinite normally ordered wedge $w$ $\in$  $F_M$ is defined similarly to the degree (\ref{e: degree}) for wedges in $ V_M^{s+nr} $. Write the vacuum vector $|M\rangle$ of  $F_M$, and the $w$ as: 
\begin{align}
& |M\rangle = u_M\wedge u_{M-1} \wedge u_{M-2} \dots \quad = z^{m_1^0}v_{\ep_1^0}\wedge z^{m_2^0}v_{\ep_2^0}\wedge \dots \quad . \\
& w = z^{m_1}v_{\ep_1}\wedge z^{m_2}v_{\ep_2}\wedge \dots \quad ; \quad m_i = m_i^0 , \ep_i = \ep_i^0 \; \text{for $ i >> 1$ }, \\  
\intertext{and define the degree $ |w | $ as: } 
& \qquad \qquad |w| = \sum_{i \geq 1} m^0_i - m_i . \label{e: sideg}
\end{align}
Let for $ k \geq 0 $ :
\begin{equation}
F_M^k = \bigoplus_{\{{\mathrm n.} {\mathrm o.}\: w \in F_M | \; |w| = k \}} \cplx w 
\end{equation}
where ``n.o.'' stands for ``normally ordered''.
The Fock space is graded with respect to this degree: $ F_M = \oplus_{k\geq 0}F_M^k $.

Define the map $ \rho_M^k : V_M^k \rightarrow F_M^k $ by:
\begin{align} 
& \text{for $ f = \{f_r\}_{r\geq 0} \in V_M^k $:} \label{e: roinf0} \\
& \rho_M^k f  = f_r \wedge | M -s-nr\rangle , \quad \text{where $ r \geq k $}.  \label{e: roinf}
\end{align}
Proposition \ref{p: prop2} (or, equivalently the eq. (\ref{e: tail})) shows that $\rho_M^k $ does not depend on the choice of $r$ in (\ref{e: roinf}) as long as $ r \geq k$. 

The following Proposition will enable us to define a level-0 action of $U_q'(\sll)$ in the Fock space.

\begin{prop}
The map $\rho_M^k $ is an isomorphism of the linear spaces $ V_M^k $ and $F_M^k$ for any $k\geq 0$.
\label{p: prop4}\end{prop}
\begin{pf}
Since $V_M^{s+nr,k}$ $\cong$ $V_M^k$ for all $ r \geq k$, it is sufficient to prove  that the map:  
\begin{equation}
V_M^{s+nk,k} \ni w^{(s+nk)} \longrightarrow w^{(s+nk)} \wedge |M - s -n k \rangle  \in F_M     \label{e: mmap}
\end{equation}
is an isomorphism of $V_M^{s+nk,k}$ and $F_M^k$.

Take a normally ordered wedge $w(\mm,\ee) \in V_M^{s+nk,k}$. This is an element of a basis of $V_M^{s+nk,k}$. Observe that the vector $w(\mm,\ee)\wedge |M - s -n k \rangle $ belongs to $ F_M^k$ and is a normally ordered wedge in $ F_M$ -- that is an element of a basis in $F_M$. This shows injectivity of the map (\ref{e: mmap}).  

Let $w: = z^{m_1}v_{\ep_1}\wedge z^{m_2}v_{\ep_2}\wedge \cdots \quad $ be a normally ordered wedge in $ F_M^k$. Applying the same reasoning as in the proof of the Proposition \ref{p: prop2} (part ii) we can show, that    
\begin{equation}
w \equiv w^{(s+nk)}\wedge |M - s -nk \rangle 
\end{equation}
where $w^{(s+nk)} \in V_M^{s+nk,k}$, and explicitely:  
\begin{equation}
w^{(s+nk)} = z^{m_1}v_{\ep_1}\wedge z^{m_2}v_{\ep_2}\wedge \cdots \wedge z^{m_{s+nk}}v_{\ep_{s+nk}}.
\end{equation}
Hence the map (\ref{e: mmap}) is surjective.
\end{pf}

Taking advantage of this Proposition we use the map $\rho_M^k $ to define in  $F_M^k $ an $U_q'(\sll)$-action $U_0$ along with an action of the commutative family $\{g_l\}$ by conjugating the actions  $U_0^{(\infty)} $ and  $g_l^{(\infty)}$ with  the isomorphism $\rho_M^k $. Obviously we have $ U_0 g_l = g_l U_0 $ for all $ l\neq 0$.
The actions in $ F_M = \oplus_{k\geq 0}F_M^k $ follow from the actions in each component $F_M^k$. This completes the definition of the level-0 action of $U_q'(\sll)$ in the Fock space $F_M$. 
\\ \mbox{} \\
 Let us summarize the results of this section. We started from $U_0^{(N)} $ -- the level-0 action of $U_q'(\sll)$ in the finite wedge product $\wedge^N V(z)$. The generators $\{x^{(N)}\}$ ( $x = E^a,F^a,K^a $ , $ a=0,\dots,n-1 $) of this action are given by (\ref{e: Kfin} - \ref{e: Ffin}) with $y_i = q^{1-N} Y^{(N)}$. Using $U_0^{(N)} $ we have defined $U_0$ -- level-0 action of $U_q'(\sll)$ in the Fock space $ F_M $ . The space $F_M$ is graded: $F_M = \oplus_{k\geq 0}F_M^k $ and by the definition the  $U_0$ preserves the degree $k$. Any vector $ w \in F_M^k $ ( $ M = s\bmod n $ , $ 0\leq s \leq n-1$ ) is represented as:  
\begin{equation}
w = f_r \wedge |M - s - nr \rangle , \label{e: wwddggee}
\end{equation}
where $ f_r \in V_M^{s+nr,k} \subset \wedge^{s+nr} V(z) $, and $ r \geq k $ . For any fixed $r \geq k$ this representation is unique by Proposition \ref{p: prop4}. The action of a generator $x \in U_q'(\sll) $ on $w$ (\ref{e: wwddggee}) is then defined as: 
\begin{equation}
x.w = (x^{(s+nr)}.f_r)\wedge |M - s - nr \rangle ,  \label{e: UactsF}
\end{equation}
and $ x.w $ does not depend on the choice of $r$ as long as $r\geq k$ again by Proposition \ref{p: prop4}. 

Similarly starting from the commutative family of operators $\{ h_l^{(N)} \}$ , $ (l =\pm 1,\pm 2,\dots ) $ defined in $\wedge^N V(z)$ we define in $F_M$ a commutative family of operators  $\{ g_l \}$, $(l =\pm 1,\pm 2,\dots )$ which also commute with the $U_0$. As in (\ref{e: UactsF}) we  prescribe the  action of $g_l$ on the wedge $w$ 
(\ref{e: wwddggee}) by:
\begin{equation}
g_l.w = (g_l^{(s+nr)}.f_r)\wedge |M - s - nr \rangle , \quad l =\pm1,\pm 2,\dots \quad ;  \label{e: gactsF}
\end{equation}
with $\{g_l^{(s+nr)}\}$ given in Proposition \ref{p: prop3} (ii). The independence this prescription on the choice of the $r$ is again due to the Proposition \ref{p: prop4}.

For computational purposes the  most convenient choice of the $r$ in (\ref{e: UactsF}, \ref{e: gactsF}) is to take it to be minimal -- that is $r=k$. We adopt this choice in the next section.

\section{Decomposition of the Fock space with respect to the level-0 action}

In this section we give the decomposition of the level-1 Fock space module of $ U_q'(\sll) $ with respect to the level-0 action  $U_0$ which was defined in sec. 4.

\subsection{Decomposition of the space $V_M^{s+nr}$}

The definition of the $ U_q'(\sll) $-action $U_0$ given in sec. 4 makes it clear, that the decomposition of the Fock space will be found once we construct the decomposition of the spaces  $V_M^{s+nr}$ and $V_M^{s+nr,k}$ with respect to the $ U_q'(\sll) $-action $U_0^{(s+nr)} $. To do this we use results of sec. 3. Let $N:=s+nr$ and $E^{\mm}$ be the subspace of $\wedge^{N} V(z)$ defined in (\ref{e: Em}). Recall that by Proposition \ref{p: finwedgedecomp} the $E^{\mm}$ is an irreducible representation of the $U_0^{(s+nr)} $. Let $\mm^0$ $\in$ $\MC_{s+nr}^n$  be the vacuum sequence (\ref{e: m0N}) associated with the integer $M$. Then the   $U_0^{(s+nr)} $-decomposition of  $V_M^{s+nr}$ is given by the following Proposition
\begin{prop}
For the $U_q'(\sll)$ modules $V_M^{s+nr}$ and $E^{\mm}$ we have
\begin{equation}
V_M^{s+nr} = \bigoplus_{\{\mm \in \MC_{s+nr}^n \; | \; m_{s+nr} \leq m^0_{s+nr}\} } E^{\mm}  \label{e: VMsdecomp}
\end{equation}\label{p: VMsdecomp}
where the $U_q'(\sll)$  action in the both sides is given by the $U^{(s+nr)}_0$.
\end{prop}
\begin{pf}
We demonstrate that the set 
\begin{equation}
B_M^{s+nr}:= \{ \phi(\mm,\ee) \equiv \Lambda(\Phi^{\mm}(\zz)\otimes {\bold v}_{\ee}) \; | \; \mm \in \MC_{s+nr}^n , \; m_{s+nr} \leq m^0_{s+nr} ; \ee \in \EC(\mm) \}  \label{e: basis}
\end{equation}  
is a base of the $V_M^{s+nr}$. 
First we note that $\phi(\mm,\ee)$ $ \in $ $B_M^{s+nr}$ implies $\phi(\mm,\ee)$ $\in $ $V_M^{s+nr}$. This follows since  the triangularity of the polynomial $\Phi^{\mm}(\zz)$ and Lemma \ref{l: lemma2} allow us to represent the $\phi(\mm,\ee)$ as 
\begin{equation}
\phi(\mm,\ee) = w(\mm,\ee)\quad + \sum_{\{\nn \in \MC_{s+nr}^n | \nn^+ \prec \mm^+ \} } \sum_{\ee' \in \EC(\nn) } c(\mm,\ee ; \nn, \ee') w(\nn,\ee').
\end{equation}
\mbox{} From the last equation it is also follows that $ B_M^{s+nr}$ is a spanning set of  $V_M^{s+nr}$. Finally the elements of  $B_M^{s+nr}$ are linearly independent by Proposition \ref{p: base}. By definition of the $E^{\mm}$ (\ref{e: Emm1}) the set $ \{ \phi(\mm,\ee) \; | \; \ee \in \EC(\mm) \} $ is a base of $ E^{\mm}$. Hence the result of the Proposition. 
\end{pf}

Since $E^{\mm}$ is homogeneous with the degree $ |\mm | $, we have 
\begin{cor}
\begin{equation}
V_M^{s+nr,k} = \bigoplus\begin{Sb} \{\mm \in \MC_{s+nr}^n \; | \; m_{s+nr} \leq m^0_{s+nr}\} \\ |\mm^0| - |\mm| = k \end{Sb} E^{\mm}  \label{e: VMskdecomp}
\end{equation}\label{c: VMskdecomp}
\end{cor}

\subsection{Decomposition of the Fock space}

Let $M = s\bmod n$, $0\leq s\leq n-1$. Write the vacuum vector $|M\rangle $ $\in $ $F_M$ as  
\begin{equation}
|M\rangle := u_{M}\wedge u_{M-1} \wedge \cdots \qquad \equiv z^{m_1^0}v_{\ep_1^0}\wedge  z^{m_2^0}v_{\ep_2^0} \wedge \cdots \qquad .
\end{equation}
The semi-infinite vacuum $\mm$-sequence associated with the $ |M\rangle $ is 
\begin{equation}
\mm^0 =(m_1^0,m_2^0,\cdots\quad ) =  (\underbrace{m^0,\dots,m^0}_s,\underbrace{m^0+1,\dots,m^0+1}_n , \underbrace{m^0+2,\dots,m^0+2}_n,\cdots \quad  ) , \quad m^0 := \frac{s-M}{n}. \label{e: m0inf} 
\end{equation}
Introduce a set $\MC^n[M] $ whose elements are  ordered semi-infinite sequences $\mm$  as follows
\begin{equation}
\MC^n[M]:= \{ \mm = (m_1 \leq m_2 \leq , \cdots \; )\; | \; \text{ $\mm$ is $n$-strict ; $m_i = m_i^0 $ for $ i>> 1$} \}.
\end{equation}

And for   $\mm \in \MC^n[M]$ define the degree $ \|\mm \| $ as 
\begin{equation}
\|\mm \| := \sum_{i \geq 1} m_i^0 - m_i \;.
\end{equation}
Note that $ \| \mm \| $ is a non-negative integer, and for a normally ordered wedge 
\begin{equation}
w = z^{m_1}v_{\ep_1}\wedge z^{m_2}v_{\ep_2}\wedge z^{m_3}v_{\ep_3}\wedge \cdots \qquad ( w \in F_M)
\end{equation}
the degree $ |w| $ (\ref{e: sideg}) is equal to the $ \|\mm \| $.

Now for $\mm \in \MC^n[M]$ denote by $\mm^{(N)}\: ( \in \MC_N^n ) $ the ordered sequence obtained from the $\mm$ by removing all except the first $N$ elements.   

Let for $\mm \in \MC^n[M]$ a linear space $ \FC_M^{\mm} \subset F_M $ be defined in the following way
\begin{equation}
\FC_M^{\mm} := E^{\mm^{(s+n \|\mm\|)}}\wedge |M-s-n\|\mm\|\rangle . \label{e: FCmm}
\end{equation}
In the last formula we have 
\begin{equation}
E^{\mm^{(s+n \|\mm\|)}} \subset V_M^{s+n \|\mm\|,\|\mm\|} \subset V_M^{s+n \|\mm\|} \subset \wedge ^{s+n \|\mm\|}V(z) ,\quad \text{and} \quad \FC_M^{\mm} \subset F_M^{\|\mm \|}.
\end{equation}

The definition of the action $U_0$ given in sec. 4  and the Proposition \ref{p: prop4} immediately lead to 
\begin{prop}
For $\mm \in \MC^n[M]$ the $\FC_M^{\mm}$ is an $U_q'(\sll)$-module with respect to the action $U_0$; the  $U_q'(\sll)$-modules  $\FC_M^{\mm}$ with the $U_q'(\sll)$-action $U_0$ and  $E^{\mm^{(s+n \|\mm\|)}}$ with the $U_q'(\sll)$-action $U_0^{(s+n\|\mm\|)}$ are isomorphic.  
\end{prop}
Note that since by Proposition \ref{p: finwedgedecomp} the $E^{\mm^{(s+n \|\mm\|)}}$ is irreducible, so is   $\FC_M^{\mm}$.

Now observing that Corollary \ref{c: VMskdecomp} and the Proposition \ref{p: prop4} imply:
\begin{equation}
\bigoplus\begin{Sb} \mm \in \MC^n[M] \\ \|\mm\| = k \end{Sb}\FC_M^{\mm}  = V_M^{(s+nk),k}\wedge |M-s - nk \rangle  = F_M^k , 
\end{equation}
we obtain for the Fock space $F_M = \oplus_{k\geq 0}F_M^k$:  
\begin{equation}
F_M = \bigoplus_{\mm \in \MC^n[M]}\FC_M^{\mm}.  \label{e: FockDec}
\end{equation}
This is the sought for decomposition of the Fock space with respect to the level-0 action $U_0$.  

Since the structure of the $E^{\mm^{(s+n \|\mm\|)}}$ as a $U_q'(\sll)$-module is known from Proposition \ref{p: finwedgedecomp} we can describe a component $\FC_M^{\mm}$ of the decomposition (\ref{e: FockDec}) by using the isomorphism of $E^{\mm^{(s+n \|\mm\|)}}$  and $\FC_M^{\mm}$. 

To do this, by analogy with (\ref{e: splitm}) for an $\mm \in \MC^n[M]$ introduce numbers $r_k$ $(k=0,1,2,\dots \;)$ by 
\begin{equation}
\mm = ( m_1 = \dots  = m_{r_1} < m_{1+r_1} =  \dots = m_{r_2} < m_{1+r_2} = \dots = m_{r_3} < \quad \cdots \quad  ), \quad \text{and $r_0:=0$}. \label{e: splitminf}
\end{equation}

As the  $U_q'(\sll)$-module the space $\FC_M^{\mm}$ then is isomorphic to the semi-infinite  tensor product: 
\begin{equation}
\FC_M^{\mm} \cong V[p^{-m_{r_1}}q^{2 r_1 - 2},r_1 ]\otimes V[p^{-m_{r_2}}q^{2 r_2 - 2},r_2 - r_1]\otimes \dots \otimes V[p^{-m_{r_k}}q^{2 r_k - 2},r_k - r_{k-1}]\otimes \quad \cdots  \label{e: Wminf}
\end{equation}
where $V[a,j]$ $(1\leq j \leq n)$ is the fundamental $U_q'(\sll)$-module defined in (\ref{e: Vaj}). Note that since $ m_i = m_i^0 $ for all but finite number of elements in $ \mm$; we have $ r_k - r_{k-1} = n $ for $ k >> 1 $, and hence the tensor product (\ref{e: Wminf}) contains only a {\em finite number} of factors different from the trivial 1-dimensional representation $V[a,n]$.  

The Drinfel'd Polynomials of $\FC_M^{\mm}$ are found from the Drinfel'd Polynomials of the representation $E^{\mm^{(s+n \|\mm\|)}}$ (see Proposition \ref{p: finwedgedecomp}), they are: 
\begin{equation}
P_i(u) =  \prod_{\{1\leq k < \infty \; | \; r_k - r_{k-1} = i \}} (u - p^{m_{r_k}}q^{-r_k-r_{k-1}}) \qquad ( i\in \{ 1,2,\dots,n-1 \}). \label{e: infwedgedecomp}
\end{equation}  
Let us remark that since $ r_k - r_{k-1} = n $ for all sufficiently large $k$, the number of factors in the product above is always finite for any $ \mm \in \MC^n[M]$.

Finally we find that the space $\FC_M^{\mm}$ is an eigenspace of the commuting Hamiltonians $g_l$ ( $ l= \pm1,\pm2,\dots $) defined by (\ref{e: gactsF}) and the  eigenvalue of $g_l$  is 
\begin{equation}
\sum_{i=1}^{\infty} ( p^{lm_i} - p^{lm_i^0})q^{2l(1-i)}. 
\end{equation}
Notice that the sum in the last expression contains only a finite number of non-zero summands due to the asymptotic condition $ m_i = m_i^0 $ for $ i >> 1$ . 

An explicit base in the space $\FC_M^{\mm}$  $( \mm \in \MC^n[M] )$ is immediately obtained from the definition (\ref{e: FCmm})  and a base of $E^{\mm^{(s+n \|\mm\|)}}$ described in (\ref{e: Emm1}).

\appendix
\section{Appendix}
In this appendix we summarize the conventions concerning the algebra $U'_{q}(\sll )$ and the Drinfel'd Polynomials adopted in this paper. In particular we discuss the proof of Lemma \ref{l: Drinfeld}.   

\subsection{Two realizations of $U'_{q}(\sll )$}

Let us recall the two realizations of $U'_{q}(\sll )$
 and the definition of the Drinfel'd Polynomials.

\begin{df} \cite{CP}
The quantum Kac--Moody algebra $U_q(g(A))$\/ associated to
a symmetric generalized Cartan matrix $A= (a_{ij})
_{i,j\in I:=\{0,1, \dots n-1\} }$\/ is the unital
associative algebra over $\cplx $ with generators $E_i$, $F_i$, $K_i^{{}\pm
1}$\/ ($i\in I$) and the following defining relations:
\begin{eqnarray}K_iK_i^{- 1} =&1& =  K_i^{ -1}K_i,\\
K_iK_j&=&K_jK_i\;,\\
K_iE_jK_i^{-1}&=& q^{ a_{ij}}E_j \; ,\\
K_iF_jK_i^{-1}&=& q^{ -a_{ij}}F_j \; ,\\
{}[E_i , F_j]&=&\delta_{ij}\frac{K_i -K_i^{-1}}{q
-q^{-1}}\;,\end{eqnarray}
\vspace{-.15in}
\begin{equation}
\sum_{r=0}^{1-a_{ij}}
 (-1)^r\left[{1-a_{ij}\atop r}\right]_{q}\;
(E_i)^rE_j (E_i)^{1-a_{ij}-r}= 0\;,\;i\ne j .
\end{equation}
\vspace{-.2in}
\begin{equation}
\sum_{r=0}^{1-a_{ij}}
 (-1)^r\left[{1-a_{ij}\atop r}\right]_{q}\;
(E_i)^rE_j (E_i)^{1-a_{ij}-r}= 0\;,\;i\ne j .
\end{equation}
\vspace{-.15in}
\begin{equation}
 \mbox{where }
[n]_{q} :=\frac{q^n -q^{-n}}{q -q^{-1}},\; \;
{}\left[{n\atop r}\right]_{q}:=\frac{[n]_{q}[n-1]_{q}\ldots {}[n-r+1]_{q}}
{[r]_{q}[r-1]_{q}\ldots [1]_{q}}.
\end{equation}
\end{df}
The coproduct $\Delta $ is given by
\begin{eqnarray}
\Delta(E_{i})& :=&  E_{i} \otimes K_{i} + 1 \otimes E_{i}, \\
\Delta(F_{i})& :=&  F_{i} \otimes 1 + K_{i}^{-1} \otimes F_{i}, \\ 
\Delta(K_{i})& :=&  K_{i} \otimes K_{i}.
\end{eqnarray}

In particular the algebra $U_{q}(\sll)$ is the algebra $U_{q}(g(A))$,
where the generalized Cartan matrix $A=(a_{ij})
_{i,j\in I}$ is
\begin{equation}
 a_{ij}= \left\{
\begin{array}{cc}
2 & (i=j) \\
-1 &(|i-j|=1 ,(i,j)=(1,n) ,(n,1)) \\
0 & (\mbox{otherwise})
\end{array} \right.
 \; \; n \geq 3 ,
\end{equation}
\vspace{-.2in}
\begin{equation}
 a_{ij}= \left\{
\begin{array}{cc}
2 & (i=j) \\
-2 &(i\neq j) 
\end{array} \right.
 \; \; n =2 .
\end{equation}
 
We put $c':=K_{0}K_{1} \dots K_{n-1}$ in $U_{q}(\sll)$,
 then $c'$ is the central in $U_{q}(\sll)$.
 We define $U'_{q}(\sll)$ as the quotient of
$U_{q}(\sll)$ by the two sided ideal generated by $c'-1$.

\begin{prop} \label{uqisom} \cite{CP}
$U'_{q}(\sll)$ is isomorphic as an algebra to the algebra $\cal{A}$\/ with
generators $E_{i,r}, \; F_{i,r}$ ($ i\in\{1,\ldots, n-1\}, r\in {\Bbb Z}$),
$H_{i,r}$ ($i\in \{1,\ldots ,n-1\}$, $r\in {\Bbb Z}\backslash \{0\}$), and
$K_i^{{}\pm 1}$, ($i\in \{1,\ldots ,n-1\}$), and the following defining
relations:
\begin{eqnarray}
K_iK_i^{- 1} =&1& =  K_i^{ -1}K_i,\\
K_iH_{j,r}&=&H_{j,r}K_i\;,\\
{}[H_{i,r}, H_{j,s}]&=& 0\;,\\
K_iE_{j,r}K_i^{-1}&=& q^{ a_{ij}}E_{j,r} \; ,\\
K_iF_{j,r}K_i^{-1}&=& q^{ -a_{ij}}F_{j,r} \; ,\\
{}[H_{i,r}, E_{j,s}]&=&{\frac1r}[ra_{ij}]_qE_{j,r+s} \;,\\
{}[H_{i,r}, F_{j,s}]&=&{-\frac1r}[ra_{ij}]_qF_{j,r+s} \;,\\
 E_{i,r+1} E_{j,s} - q^{ a_{ij}}E_{j,s} E_{i,r+1}&=&q^{ a_{ij}}
E_{i,r}E_{j,s+1} - E_{j,s+1}E_{i,r},\\
 F_{i,r+1} F_{j,s} - q^{- a_{ij}}F_{j,s} F_{i,r+1}&=&q^{- a_{ij}}
F_{i,r}F_{j,s+1} - F_{j,s+1}F_{i,r},\\
{}[E_{i,r} , F_{j,s}]&=&\delta_{ij}\frac{\Phi_{i,r+s}^+ -\Phi_{i,r+s}^-
}{q -q^{-1}}\;, \label{phi+} \end{eqnarray}
\vspace{-.25in}
\begin{equation}
\sum_{\pi\in S_p}\sum_{k=0}^{p}
 (-1)^k\left[{p\atop k}\right]_{q}\; E_{i,r_{\pi(1)}}\ldots
E_{i,r_{\pi(k)}}E_{j,s}E_{i,r_{\pi(k+1)}}\ldots
E_{i,r_{\pi(p)}}= 0\;,\;i\ne j ,
\end{equation}
\vspace{-.25in}
\begin{equation}
\sum_{\pi\in S_p}\sum_{k=0}^{p}
 (-1)^k\left[{p\atop k}\right]_{q}\; F_{i,r_{\pi(1)}}\ldots
F_{i,r_{\pi(k)}}F_{j,s}F_{i,r_{\pi(k+1)}}\ldots
F_{i,r_{\pi(p)}}= 0\;,\;i\ne j ,
\end{equation}
for all sequences $(r_1,\ldots,r_p)\in{\Bbb Z}^p$,
where $p = 1-a_{ij}$\/ and the elements $\Phi_{i,r}^{{}\pm{}}$\/ are determined
by equating coefficients of powers of $u$\/ in the formal power series
\begin{equation}
\Phi_{i}^{\pm }(u):= 
\sum_{r=0}^{\infty}\Phi_{i, {{}\pm r}}^{{}\pm{}}u^{{}\pm r} = K_i^{{}\pm
1}exp({}\pm (q-q^{-1})\sum_{s=1}^{\infty} H_{i, {}\pm s} u^{{}\pm s}) .
\end{equation}
The generators of $\cal{A}$ are called Drinfel'd generators.
\end{prop}

The isomorphisms $\tilde{f} :U'_{q}(\sll)
 \rightarrow \cal{A}$ are not determined uniquely.
In this paper we fix one isomorphism
 $f :U'_{q}(\sll) \rightarrow \cal{A}$
\begin{equation}
f(E_i) = E_{i,0}, \; \; 
f(F_i) = F_{i,0},
\;\; f(K_i^{{}\pm 1}) = K_i^{{}\pm 1},
\end{equation}
for $i\in\{1,\ldots, n-1\}$, and
\begin{eqnarray} & f(K_0^{{}\pm 1}) = (K_1K_2\ldots K_{n-1})^{{}\mp 1}, & \\
& f(E_0)=(-1)^{m-1}q^{-(n-2)/2}[F_{n-1,0},[F_{n-2, 0},\ldots 
& \label{e0} \\
&  \ldots ,[F_{m+1,0} ,[F_{1,0}
,\ldots ,[F_{m-1,0} ,F_{m,1} ]_{q^{1/2}}\ldots ]_{q^{1/2}}f(K_0), &
 \nonumber \\
& f(F_0)= \mu f(K_0^{-1})[E_{n-1,0},[E_{n-2, 0},\ldots & \\ 
& \ldots ,[E_{m+1,0},[E_{1,0}
,\ldots ,[E_{m-1,0},E_{m,-1}]_{q^{1/2}}\ldots ]_{q^{1/2}},& \nonumber
\end{eqnarray}
where $\mu\in\cplx ^{\times } $\/ is determined by
\begin{equation}
[f(E_0), f(F_0)] = \frac{f(K_0) - f(K_0^{-1})}{q-q^{-1}},
\end{equation}
and $\mbox{} [ a,b]_{q^{1/2}} := q^{1/2} ab -q^{-1/2}ba$. \\

\subsection{Drinfel'd Polynomials}

 Now we will  calculate the Drinfel'd Polynomials. 

Let $W$ be the representation of $U'_{q}(\sll)$.
 $W$ is said to be of type 1 if
\begin{equation}
W=\bigoplus_{{\bf\mu}\in {\Bbb Z}^{n}}W_{{\bf\mu}},
\end{equation}
where $W_{\mu} =\{w\in W|k_i.w = q^{\mu(i)} w\}$. 
\begin{prop} 
\cite{CP}
Let $W$\/ be a finite-dimensional irreducible $U'_{q}(\sll)$--module of type 1.
Then, 

(a) $W$\/ is generated by a vector $w_0$\/ satisfying
\begin{equation}
E_{i,r}. w_0 = 0, \;\; \Phi_{i,r}^{{}\pm{}}.w_0 =\phi_{i,r}^{{}\pm{}}w_0 
\end{equation}
for all $i\in\{1,\ldots ,n\}$, $r\in{\Bbb Z}$, and some
$\phi_{i,r}^{{}\pm{}}\in\cplx $.

(b) There exist unique monic polynomials $P_1(u),\ldots, P_{n-1}(u)$\/
(depending on $W$)  such that the $\phi_{i,r}^{{}\pm{}}$\/ satisfy
\begin{equation}
\sum_{r=0}^{\infty}\phi_{i,r}^+u^r =q^{{\rm
deg}\:P_i}\frac{P_i(q^{-2}u)}{P_i(u)} =\sum_{r=0}^{\infty}\phi_{i,r}^-u^{-r},
\end{equation}
in the sense that the left and right-hand sides  are the Laurent expansions of
the middle term about $0$\/ and $\infty$\/ respectively.
Assigning to $W$\/ the corresponding $n-1$--tuple of polynomials defines a one to
one correspondence between the isomorphism classes of finite--dimensional
irreducible $U'_{q}(\sll)$--modules of type 1 and the set of $n-1$--tuples of monic
polynomials in one variable $u$. We define the polynomials
 $P_1(u),\ldots, P_{n-1}(u)$ to be the Drinfel'd Polynomials.
\end{prop} 

\begin{rmk}
If we change the isomorphism of the Proposition \ref{uqisom},
the Drinfel'd Polynomials may be changed.
\end{rmk}

As a consequence of this Proposition, we get

\begin{cor}
\label{dp2} \cite{CP}
Let $W$\/ be a finite-dimensional irreducible representation of  $U'_{q}(\sll)$
with associated polynomials $P_i$. Set $\lambda =({\rm deg}\;P_1,\ldots, {\rm
deg}\; P_n)$. Then $W$\/ contains the irreducible
 $U_{q}(\sln)$--module  $V(\lambda)$\/
with multiplicity one. Further, if $V(\mu)$\/ is any other
 $U_{q}(\sln)$--module
occurring in $W$, then $\lambda \succeq \mu$. 
\end{cor}

Let $\Lambda _{j}$ be the $j$-th fundamental weight of $\sln$, 
 then $V(\Lambda _{j})$ is an irreducible representation of $U_{q}(\sln)$.
If $V(\Lambda _{j})$ is also the representation of
 $U'_{q}(\sll) \; ( \supset U_{q}(\sln) )$, then $V(\Lambda _{j})$ 
is irreducible as a $U'_{q}(\sll)$-module and by the Corollary \ref{dp2},
 and the Drinfel'd Polynomials of $V(\Lambda _{j})$ is 
\begin{equation}
 P_i(u) =\left\{\begin{array}{ll} u-\tilde{a}\; \;& {\rm if}\; i =j,\\
1\;\; & {\rm otherwise,} \end{array}\right. \label{repofm}
\end{equation}
for some constant $\tilde{a}$. We define the representation determined by
 (\ref{repofm}) as $V(\Lambda _{j}; \tilde{a})$.

We need the following Lemma.
\begin{lemma}
\label{e0lem} \cite{CP}
 Let $v_{\Lambda_j}$\/ be the $U_{q}(\sln)$--highest weight vector
 in $V(\Lambda_j ;\tilde{a})$,  where $m\in\{1,\ldots ,n-1\}$, $\tilde{a}
\in\cplx^{\times }$. Then,
\begin{equation}
E_0.v_{\Lambda_j} = (-1)^{j-1}q^{-1}\tilde{a}^{-1}F_{n-1}F_{n-2}\ldots
F_{j+1}F_1\ldots F_j. v_{\Lambda_j}.
\end{equation}
\end{lemma}

\begin{rmk}
This Lemma and the Lemma 6.4 in the paper \cite{CP} are different
 because the isomorphism between the realization of the Chevalley generators  
and the realization of the Drinfel'd generators are different.
\end{rmk} \\
{\it{proof}} \hspace{.3in}
 Using the fact that the weight spaces of $V(\Lambda_j ,\tilde{a})$\/ 
as a $U_{q}(\sln)$--module are all one--dimensional and 
$K_{i}F_{j,1}K_{i}^{-1} = q^{- a_{ij}} F_{j,1}$, we get
\begin{equation}
F_{j,1} .v_{\Lambda_j} = b F_j .v_{\Lambda_j}
\end{equation}
for some $b\in \cplx$. \mbox{} From the relation (\ref{phi+}), we get
\begin{equation}
\Phi_{j,1}^+.v_{\Lambda_j} = b(q-q^{-1})v_{\Lambda_j}.
\end{equation}
Hence, from the definition of the Drinfel'd Polynomials, we have
\begin{equation} 
q(q^{-2}u-\tilde{a}) = (u-\tilde{a})(q+b(q-q^{-1})u +O(u^2)),
\end{equation}
so that $b =\tilde{a}^{-1}$.
 Finally, from the relation (\ref{e0}), we find that
\begin{equation} 
E_0.v_{\Lambda_j} = (-1)^{j-1} q^{-1}\tilde{a}^{-1}F_{n-1}F_{n-2}\ldots
F_{j-1}F_1\ldots F_j. v_{\Lambda_j}.
\end{equation}
\hfill\halmos

 \mbox{} From now on, we calculate the Drinfel'd Polynomials of $V[a,j]$
 (\ref{e: Vaj}).
$V[a,j]$ is the highest weight representation as a $U_{q}(\sln )$-module,
 and the highest weight is $\Lambda_{j} $ (the $j$-th fundamental weight).
First we can check that $ {} [ v_{1} \otimes v_{2} \otimes \dots \otimes
v_{j} ] $ is the highest weight vector of $V[a,j]$
 as a $U_{q}(\sln )$--module. Because of the argument before
 the Lemma \ref{e0lem}, the Drinfel'd Polynomials of $V[a,j]$ are 
\begin{equation}
 P_i(u) =\left\{\begin{array}{ll} u-\tilde{a}\; \;& {\rm if}\; i =j,\\
1\;\; & {\rm otherwise,} \end{array}\right. 
\end{equation}
for some constant $\tilde{a}$. To determine $\tilde{a}$,
 we observe how the Chevalley generators act.
\begin{eqnarray}   
& E_{0} \cdot ( [ v_{1} \otimes v_{2} \otimes \dots \otimes v_{j} ] ) & \\
& = [ E_{0} v_{1} \otimes K_{0} v_{2} \otimes \dots \otimes K_{0} v_{j} ]+[
v_{1} \otimes E_{0} v_{2} \otimes K_{0} v_{3} \otimes \dots \otimes K_{0} v_{j}
] +  \dots  & \nonumber \\
& = a [ v_{n} \otimes v_{2} \otimes \dots \otimes v_{j} ] . & \nonumber
\end{eqnarray}
\vspace{-.15in}
\begin{eqnarray}
&  F_{n-1} \dots F_{j+1} F_{1} \dots F_{j} \cdot
 ( [ v_{1} \otimes v_{2} \otimes \dots \otimes v_{j} ] ) & \\
&  =  F_{n-1} \dots F_{j+1} F_{1} \dots F_{j-1}
 \cdot (\dots +[ K_{j}^{-1} v_{1} \otimes \dots \otimes K_{j}^{-1} v_{j-1} 
\otimes F_{j} v_{j} ] ) \nonumber & \\
&  = \dots =  F_{n-1} \dots F_{j+1} \cdot
 ( [ v_{2} \otimes v_{3} \otimes \dots \otimes v_{j+1} ]) & \nonumber \\
&  =  F_{n-1} \dots F_{j+2} \cdot 
 ( \dots +[ K_{j+1}^{-1}v_{2} \otimes K_{j+1}^{-1} v_{3}
 \otimes \dots \otimes F_{j+1} v_{j+1} ] )
&  \nonumber \\ 
&  = \dots = ([v_{2} \otimes \dots \otimes v_{j}\otimes v_{n} ]). & \nonumber
\end{eqnarray}

On the other hand we can check 
\begin{equation} 
v_{k} \otimes v_{l} + q v_{l} \otimes v_{k} \in Ker(q+S^{-1}) \;
\mbox { for } k<l .
\end{equation}
Then 
\begin{eqnarray}
&  \mbox{} [v_{2} \otimes \dots \otimes v_{j} \otimes v_{n}] 
 = (-q) [v_{2} \otimes \dots v_{j-1} \otimes v_{n} \otimes v_{j}] & \\
& = (-q)^{2} [v_{2} \otimes \dots \otimes 
v_{n} \otimes v_{j-1} \otimes v_{j}] =
\dots = (-q)^{j-1} [ v_{n} \otimes v_{2} \otimes \dots \otimes v_{j} ] . & 
\nonumber
\end{eqnarray}
\mbox{} From the Lemma \ref{e0lem}, we get $\tilde{a} = q^{j-2} a^{-1}$.
So we have proved the Lemma \ref{l: Drinfeld}.

Next we verify the Proposition \ref{p: finwedgedecomp}.
 Basically we use the following Proposition.  
\begin{prop} \cite{CP}
If $V , \; W , \; V \otimes W$ are all irreducible as a $U'_{q}(\sll)$--module,
 whose Drinfel'd Polynomials are $P_{V,i} (u), \; P_{W,i} (u), \;
P_{V \otimes W ,i} (u) \; \; (i\in \{ 1,2, \dots ,n-1 \}) $, then 
\begin{equation}
P_{V \otimes W ,i} (u)= P_{V,i} (u) \cdot P_{W,i} (u).
\end{equation}
\end{prop}

Using the fact written in the paper \cite{AK}, we can check that
 our representation $W^{\mm } $ (\ref{e: Wm1}) 
 is irreducible if $\alpha \in \real \setminus {\Bbb Q}_{\leq 0} \; \;
(p=q^{-2\alpha })$. We explain the reason briefly.
Let $\Lambda _{k_{i}}$ be the $k_{i}$-th fundamental weight of $\sln $.
The key propositions are \\
$\circ \; V_{1}(\Lambda _{k_{1}} ; a_{1}) \otimes 
V_{2}(\Lambda _{k_{2}} ; a_{2})$
is irreducible if and only if the intertwiner $R( V_{1}, V_{2}):
V_{1}(\Lambda _{k_{1}} ; a_{1}) \otimes
 V_{2}(\Lambda _{k_{2}} ; a_{2}) \rightarrow
V_{2}(\Lambda _{k_{2}} ; a_{2}) \otimes 
V_{1}(\Lambda _{k_{1}} ; a_{1})$ and $R(V_{2} , V_{1})$ has no pole. \\
$\circ \; V_{1}(\Lambda _{k_{1}} ; a_{1}) \otimes \dots \otimes
 V_{l}(\Lambda _{k_{l}} ; a_{l})$ is irreducible if and only if 
$ V_{i}(\Lambda _{k_{i}} ; a_{i}) \otimes 
V_{j}(\Lambda _{k_{j}} ; a_{j}) $ is irreducible for all $i<j$.

\end{document}